\documentclass[reqno, 12pt, oneside]{amsart}

\usepackage{amsmath,amsthm,amsfonts,amssymb, latexsym, eepic, xypic, fancyhdr,cancel, geometry}
\usepackage[dvips]{epsfig}
\usepackage[usenames]{color}
\usepackage{enumerate}

\newtheorem{theorem}{Theorem} 
\newtheorem{lemma}[theorem]{Lemma}
\newtheorem{definition}[theorem]{Definition}
\newtheorem{corollary}[theorem]{Corollary}
\newtheorem{proposition}[theorem]{Proposition}    
\newtheorem{remark}[theorem]{Remark}

\renewcommand{\L}{\mathcal{L}}
\renewcommand{\div}{\operatorname{div}}
\renewcommand {\H} {\mathcal{H}}
\newcommand{\C}{\operatorname{\mathcal{C}}}
\newcommand {\R} {\operatorname {\mathbb{R}}}

\newcommand{\Ric}{\operatorname{Rc}}
\newcommand{\supp}{\operatorname{supp}}
\newcommand {\graph} {\operatorname{graph}}
\newcommand {\tr} {\operatorname{tr}}

\begin{document}

\title[The Jang Equation Reduction]{The Jang equation reduction of the spacetime positive energy theorem in dimensions less than eight}

\author[Eichmair]{Michael Eichmair}
\address{Michael Eichmair, ETH Z\"urich, Departement Mathematik, 8092 Z\"urich, Switzerland}
\email{michael.eichmair@math.ethz.ch}

\thanks{The author gratefully acknowledges the support of NSF grant DMS-0906038 and of SNF grant 2-77348-12.}

\date{\today}

\begin{abstract} We extend the Jang equation proof of the positive energy theorem due to R. Schoen and S.-T. Yau \cite{Schoen-Yau:1981-pmt2} from dimension $n=3$ to dimensions $3 \leq n <8$. This requires us to address several technical difficulties that are not present when $n=3$. The regularity and decay assumptions for the initial data sets to which our argument applies are weaker than those in \cite{Schoen-Yau:1981-pmt2}. 
\end{abstract}
\maketitle


\section{Introduction}

We begin by recalling the relevant definitions. 

\begin{definition} [Cf. \protect{\cite[Definitions 3 and 4]{spacetimePMT}}]  \label{def:AF} Let $n \geq 3$. An initial data set $(M, g, k)$ consists of a complete boundaryless $n$-dimensional $\C^{2}$ Riemannian manifold $(M, g)$ and a $\C^{1}$ symmetric  $(0, 2)$-tensor $k$. The local mass density $\mu$ and the local current density $J$ of an initial data set $(M, g, k)$ are defined as 
\begin{eqnarray*}
\mu := \frac{1}{2} \left( R_g - |k|_g^2 + \tr_g(k)^2 \right) \text{ and } J := \div_g \left( k - \tr_g(k) g \right). 
\end{eqnarray*} 
Here, $R_g$ denotes the scalar curvature of $(M, g)$. 
Given $\ell \in \{2, 3, \ldots\}$, $p>n$, $q \in (\frac{n-2}{2}, n-2)$, $\alpha \in (0, 1- \frac{n}{p}]$, and $q_0 >0$, we say that an initial data set $(M, g, k)$ is asymptotically flat of type $(\ell, p, q, q_0, \alpha)$, if $g \in \C^{\ell, \alpha} (M)$, if $k \in \C^{\ell - 1, \alpha} (M)$, and if there exists a compact set $K \subset M$,  a closed coordinate ball $\bar B \subset \R^n$, and a $\C^{\ell+1, \alpha}$ diffeomorphism $x = (x_1, \ldots, x_n): M \setminus K \to \R^n \setminus \bar B$ such that 
\begin{eqnarray}
\label{eqn:decaymetric} g_{ij} - \delta_{ij} \in W^{\ell, p}_{-q} (\R^n \setminus \bar B), \ k_{ij} \in W^{\ell-1, p}_{-q-1} (\R^n \setminus \bar B), \\ \text{ and } (\mu, J) = O^{\ell-2, \alpha} (|x|^{-n-q_0}) \text{ as } |x| := \sqrt{x_1^2 + \ldots + x_n^2} \to \infty. \nonumber
\end{eqnarray}
An initial data set $(M, g, k)$ is said to have harmonic asymptotics of type $(\ell, \alpha)$ if the diffeomorphism $x = (x_1, \ldots, x_n)$ above can be chosen such that in addition there exists a positive function $\varphi \in \C^{\ell, \alpha} (M)$ and a vector field $Y \in \C^{\ell, \alpha} (\Gamma (T M))$ for which 
\begin{eqnarray} \label{eqn:harmonicallyflat}
g_{ij} = \varphi^{\frac{4}{n-2}} \delta_{ij} \text { and } k_{ij} = \varphi^{\frac{2}{n-2}} \left( Y_{i, j} + Y_{j, i} \right)  \text{ on } \R^n \setminus \bar B, 
\end{eqnarray}
and such that there exist constants $a, b_1, \ldots, b_n$ for which 
\begin{eqnarray} \label{eqn:uY}
\varphi(x) - 1 - a |x|^{2-n} &=& O^{\ell, \alpha} (|x|^{1-n}) \text{ and } \\ Y_i - b_i |x|^{2-n} &=& O^{\ell, \alpha} (|x|^{1-n}) \text{ for all } i = 1, \ldots, n. \nonumber
\end{eqnarray}
\end{definition}

\noindent The precise definitions of the (standard) weighted Sobolev and H\"older spaces in the preceding definition are given in \cite[Section 2]{spacetimePMT}. Note that for an initial data set that has harmonic asymptotics of type $(\ell, \alpha)$ we have that $(\mu, J) = O^{\ell-2, \alpha} (|x|^{-n-1})$ and  $\tr_g (k) = O^{\ell, \alpha}(|x|^{-n})$. For convenience and definiteness, we agree that $B = B (0, \frac{1}{3})$, and we denote by $|\cdot|$ an extension of the Euclidean coordinate radius function on $x^{-1} (B (0, \frac{1}{2})) \subset M$ to a function in $\C^{\ell + 1, \alpha} (M)$ that is bounded below by $\frac{1}{3}$. \\

We recall two facts about asymptotically flat initial data sets $(M, g, k)$. First, the ADM energy $E$ and the ADM momentum tensor $P$ of $(M, g, k)$ are well-defined by the expressions
\begin{eqnarray}
\label{defn:energy} E = \frac{1}{2(n-1) |\mathbb{S}^{n-1}|} \lim_{r \to \infty} \int_{|x| = r} \sum_{i, j=1}^n \left(\partial_j g_{ij} - \partial_j g_{ii}\right) \frac{x_j}{r}  d \H^{n-1}_\delta \text{ and } \\
P_i = \frac{1}{(n-1) |\mathbb{S}^{n-1}|} \lim_{r \to \infty} \int_{|x| = r} \sum_{j=1}^n \left( k_{ij} - \tr_g (k) g_{ij} \right) \frac{x_j}{|x|} d \H^{n-1}_\delta \\ \text{ for } i = 1, \ldots, n. \nonumber
\end{eqnarray} 
We refer the reader to \cite[Section 5.1]{ADM:1962, ADM:2004} by R. Arnowitt, S. Deser, and C. W. Misner for the original definition of the energy and (linear) momentum of an initial data set and the paper \cite{Bartnik:1986} by R. Bartnik concerning the independence of these quantities of the asymptotic coordinate system in which they are computed. When the initial data set $(M, g, k)$ has harmonic asymptotics, then, cf. \cite{spacetimePMT},
\begin{eqnarray} \label{eqn:EPab}
E = \frac{(n-2)}{2} a \ \ \ \text{ and } \ \ \ P_i = - \frac{n-2}{n-1} b_i \ \ \  \text{ for }  i = 1, \ldots, n, 
\end{eqnarray}
where $a$ and $b_1, \ldots, b_n$ are as in (\ref{eqn:uY}). 
 
Second, we will use the following density theorem.

\begin{theorem} [\protect{\cite[Theorem 18]{spacetimePMT}}] \label{thm:density} Let $(M, g, k)$ be an asymptotically flat initial data set of type $(\ell, p, q, q_0, \alpha)$ with asymptotic coordinate chart $x : M \setminus K \to \R^n \setminus \bar B$. Assume that the dominant energy condition $\mu \geq |J|_g$ holds on $M$. There exist compact sets $K_j \subset M$ each containing $K$ and a sequence of asymptotically flat initial data $(g^j, k^j)$ that have harmonic asymptotics of type $(\ell, \alpha)$ on $M \setminus K_j$ in the chart $x$ such that $||(g, k)-(g^j, k^j)||_{W^{\ell, p}_{-q} \times W^{\ell-1, p}_{-1-q}} \to 0$, $||(\mu, J) - (\mu^j, J^j)||_{\C^{\ell-2, \alpha}_{-n-q_0}} \to 0$, and $(E^j, P^j_i) \to (E, P_i)$ as $j \to \infty$, and such that the strict dominant energy condition $\mu^j> | J^j|_{g^j}$ holds for all $j = 1, 2, \ldots$. Given $q_0' \geq q_0$ we can choose $(g^j, k^j)$ such that $(\mu^j, J^j) = O^{\ell-2, \alpha}(|x|^{-n-q_0'})$. 
\end{theorem}

The purpose of this paper is to prove the following theorem: 

\begin{theorem} [Spacetime positive energy theorem] \label{thm:main} Let $3 \leq n < 8$ and let $(M, g, k)$ be an asymptotically flat initial data set of type $(2, p, q, q_0, \alpha)$ of dimension $n$ such that the dominant energy condition $\mu \geq |J|_g$ holds on $M$. Then $E \geq 0$. 
If $E = 0$ and if
\begin{eqnarray} \label{eqn:trace3} 
\tr_g (k) = g^{ij} k_{ij} = O(|x|^{-\gamma}) \text{ for some } \gamma >2,
\end{eqnarray}
then $(M, g, k)$ is Cauchy initial data for Minkowski space, i.e. $(M, g)$ can be embedded isometrically as a spacelike hypersurface into $(\R^n \times \R, dx_1^2 + \ldots + dx_n^2 - dx_{n+1}^2)$ with second fundamental form $k$. 
\end{theorem}

A few remarks are in order.  First, note that condition (\ref{eqn:trace3}) is redundant when $n>3$. The special case  of Theorem \ref{thm:main} when $k \equiv 0$ was proven by R. Schoen and S.-T. Yau in \cite{Schoen-Yau:1979-pmt1} in dimension $n=3$, and then in dimensions $3\leq n<8$ in \cite{Schoen-Yau:1979-pat}. This is known as the ``Riemannian" or ``time-symmetric" positive energy theorem. They established the special case of Theorem \ref{thm:main} where $n=3$ in \cite{Schoen-Yau:1981-pmt2}. Their proof is  by  reduction to the Riemannian positive energy theorem using the Jang equation. In joint work \cite{spacetimePMT} with L.-H. Huang, D. Lee, and R. Schoen, the author has recently proved the stronger conclusion $E \geq |P|$, i.e. the full positive mass theorem, when $3 \leq n < 8$. Under the additional topological assumption that $M$ is spin, this had been known previously in all dimensions $n \geq 3$ by work of E. Witten \cite{Witten:1981}, see also \cite{Bartnik:1986} for details. The introduction of \cite{Witten:1981} contains a detailed discussion of earlier results in the physics literature, notably the paper \cite{Deser-Teitelboim:1977} by S. Deser and C. Teitelboim and the paper \cite{Grisaru:1978} by M.T. Grisaru; see also \cite{Deser:1983}. \\

We pointed out in \cite{spacetimePMT} that the full positive mass theorem, $E \geq |P|$, can be obtained from Theorem \ref{thm:main} via a reduction argument based on our density theorem \cite[Theorem 18]{spacetimePMT} and the boost argument of D. Christodoulou and N. \'O Murchadha \cite{CO:1981}. \\

The extension of the proof in \cite{Schoen-Yau:1981-pmt2} to higher dimensions is not all straightforward. The purpose of the present paper is to identify and resolve the ensuing difficulties. From this point of view, the reader should look at this work as a compendium to \cite{Schoen-Yau:1981-pmt2}.\\

An obvious and necessary difference with \cite{Schoen-Yau:1981-pmt2} is that instead of stability based curvature estimates (which are available only when $n \leq 5$) we rely on techniques from geometric measure theory that have been introduced by the author to the study of the Jang equation and the (related) study of marginally outer trapped surfaces in \cite{Eichmair:2009-Plateau, Eichmair:2010}. A second difference with \cite{Schoen-Yau:1981-pmt2} is the analysis of the conformal structure of certain complete exterior solutions of the Jang equation that have cylindrical blow ups or blow downs along the boundary components of their domains. If the strict dominant energy condition holds on $(M, g, k)$ and $n=3$, the topology of these boundary components is always that of $\mathbb{S}^2$ and the cylindrical ends can be compactified conformally so they become equivalent to the flat metric on a punctured ball. When $n>3$, all we know about these boundary components  is that they have positive Yamabe type. In Propositions \ref{prop:conformalfactor} and \ref{prop:blowupfunction} we deal with this difference from the case $n=3$ in the conformal analysis. Since all these points are technical and delicate in nature, we supply additional details to many of the steps in \cite{Schoen-Yau:1981-pmt2} throughout the paper. Finally, our regularity and decay assumptions on the initial data $(M, g, k)$ in Theorem \ref{thm:main} are weaker than those of \cite{Schoen-Yau:1981-pmt2}. For example, in order to characterize the equality case in Theorem \ref{thm:main}, the metric is assumed to be in $\C^4$ with appropriate decay in \cite{Schoen-Yau:1981-pmt2}. \\

We point out that Theorem \ref{thm:main} also applies to asymptotically flat initial data sets with more than one end. In this case, the diffeomorphism $x : M \setminus K \to \coprod_{i=1}^N (\R^n \setminus B)_i$ in Definition \ref{def:AF} is onto the disjoint union of $N$ exterior regions $\{(\R^n \setminus B)_i\}_{i=1}^N$. The assertion is that the ADM-energy of each end is non-negative, and that $(M, g, k)$ embeds as a spacelike hypersurface into Minkowski space (and hence in particular has the topology of $\R^n$) if the ADM-energy of any one end vanishes. In \cite{Schoen-Yau:1981-pmt2}, this is shown when $n=3$ by fixing one of the ends and carrying the remaining ends along throughout the proof. This requires several additional technical steps. We point out an alternative route here. A chosen end can be separated from all the other ends of $(M, g, k)$ by an outermost marginally trapped surface, see \cite{Andersson-Metzger:2009} for the case $n=3$, and \cite{Eichmair:2010} for a different proof that works for all $3 \leq n <8$. The main result in \cite{Metzger:2010} (see alternatively \cite [Remark 4.1] {Eichmair:2010}) shows that one can find a complete connected hypersurface $\Sigma = \graph(f_\Sigma, U_\Sigma)$ as in (\ref{exteriorsolution}) of Proposition \ref{prop:limitanalysis} such that $U_\Sigma$ contains the exterior region of the chosen end, and such that $U_\Sigma$ lies beyond the outermost marginally trapped surface that separates the chosen end from the other ends. The rest of the proof can then proceed exactly as in the case of one end. 

\subsection*{Acknowledgements} This paper is based on my thesis at Stanford University, and I would like to thank my advisor Richard Schoen and also Simon Brendle, Leon Simon, and Brian White for their unfailing support. My work -- and I -- owe a great deal to their excellent guidance and teaching. The perturbation argument to the strict dominant energy condition in my thesis is flawed. Here we substitute a  density argument from my recent paper with Lan-Hsuang Huang, Dan Lee, and Richard Schoen for this step. The exposition here also benefits from my collaboration with Lars Andersson and Jan Metzger on our survey article \cite{Andersson-Eichmair-Metzger:2010}. I would like to extend my special thanks to Jan Metzger and Anna Sakovich for proofreading a preliminary version of this manuscript carefully. I am also grateful to Robert Beig, Piotr Chrus\'ciel, Justin Corvino, Mattias Dahl, Greg Galloway, and Dan Pollack for their interest, discussions, and encouragement, and to Stanley Deser for his interest and comments. Finally, I would like to extend my appreciation and thanks to the referees. Their constructive criticism has led to a significant improvement of this paper.    


\section{Geometric solutions of the Jang equation} \label{sec:geometricsolutions}

Let $U \subset M$ be an open subset and let $u \in \C^2_{loc} (U)$. We let
\begin{eqnarray*}
H(u) := \div_g  \frac{\nabla_g u}{\sqrt{1 + |d u|_g^2}}  = \left( g^{ij} - \frac{\nabla_g^i u \nabla_g^j u}{1 + |d u|^2_g} \right) \frac{(\nabla^2_{g} u)_{ij} }{\sqrt{1 + |d u|_g^2}}
\end{eqnarray*}
and 
\begin{eqnarray*}
\tr_g(k)(u) := \tr_g(k) -  (1 + |d u|^2_g)^{-1}k(\nabla_g u, \nabla_g u) = \left( g^{ij} - \frac{\nabla^i_g u  \nabla^j_g u}{1 + |d u|^2_g} \right) k_{ij}. 
\end{eqnarray*}
We recall that $H(u)$ is the scalar mean curvature of $\graph(u, U) \subset (M \times \R, g + dt^2)$ computed as the tangential divergence of the downward pointing unit normal, and that $\tr_g(k)(u)$ is the trace over the tangent space of $\graph(u, U)$ of the tensor $k$ extended trivially in the vertical direction to $M \times \R$. \\

In this section we study the existence of solutions $f_\tau  \in \C^{2, \alpha}_{loc}(M)$ of the prescribed mean curvature equations
\begin{eqnarray} \label{eqn:Jangequation}
H(f_\tau) - \tr_g(k)(f_\tau) = \tau f_\tau
\end{eqnarray}
for $\tau >0$ and geometric limits of graphs of solutions $\graph(f_\tau, M) \subset M \times \R$ as $\tau \searrow 0$.\\

\noindent Equation (\ref{eqn:Jangequation}) with $\tau = 0$ was first considered by P.-S. Jang in \cite{Jang:1978} in an attempt to prove the spacetime positive energy theorem using a geometric flow argument. We will refer to it as the Jang equation. The assumption that an entire solution of the Jang equation with decay in the asymptotically flat end of $M$ exists underlies the strategy of \cite{Jang:1978}. The capillarity regularization (\ref{eqn:Jangequation}) of the Jang equation was introduced by R. Schoen and S.-T. Yau in \cite{Schoen-Yau:1981-pmt2}. They discovered that there are natural obstructions, namely closed marginally trapped surfaces in the initial data set $(M, g, k)$, to the existence of entire solutions of the Jang equation, while solutions of (\ref{eqn:Jangequation}) always exist. R. Schoen and S.-T. Yau introduced  minimal surface techniques to study geometric limits of the submanifolds $\graph(f_\tau, M) \subset M \times \R$ as $\tau \searrow 0$ when $n=3$. Subsequential  limits exist, and they contain complete graphical components that are solutions of the Jang equation on exterior regions of $M$, with blow up or blow down along the boundary components of their domains. \\

A basic ingredient in \cite{Schoen-Yau:1981-pmt2} are curvature estimates for $\graph(f_\tau, M)$ analogous to those for stable minimal hypersurfaces in \cite{Schoen-Simon-Yau:1975}. These estimates are available so long as $n \leq 5$. In \cite{Eichmair:2009-Plateau} and \cite{Eichmair:2010} the author observed that both the stability based regularity theory of R. Schoen and L. Simon \cite{Schoen-Simon:1981} and the theory of almost minimizing boundaries as introduced and studied in \cite{Almgren:1976, Tamanini:1982, Duzaar-Steffen:1993} are available to extend the geometric theory of the Jang equation and, more importantly, that of marginally trapped surfaces to all dimensions $n \geq 3$. We refer the reader to the author's survey article with L. Andersson and J. Metzger \cite{Andersson-Eichmair-Metzger:2010} for a discussion of related developments.  \\

The following proposition is a simple extension of the argument in \cite[p. 248]{Schoen-Yau:1981-pmt2} to dimensions $n \geq 3$: 

\begin{proposition} [Cf. p. 248 in \protect{\cite{Schoen-Yau:1981-pmt2}} when $n=3$]\label{prop:barriers} Fix $\beta \in (2, n)$. For $\Lambda \geq 1$ define $${}_\Lambda b (r) := \Lambda \int_{\frac{r}{\Lambda}}^\infty \frac{ds}{\sqrt{s^{2(\beta -1)} - 1}} \text{ on } [\Lambda, \infty).$$ Then ${}_\Lambda b$ is positive, continuous, smooth on $(\Lambda, \infty)$, and $\frac{d {}_\Lambda b}{dr} (r)$ tends to $-\infty$ as $r \searrow \Lambda$. There exists a constant $c = c(\beta) \geq 1$ such that ${}_\Lambda b (r) \leq c \Lambda (r/\Lambda)^{2- \beta}$ and such that ${}_\Lambda b (\Lambda) \geq c^{-1} \Lambda$. 

Let $(M, g, k)$ be an asymptotically flat initial data set of dimension $n$, $n \geq 3$, with $g^{ij} k_{ij} = O(|x|^{-\beta})$. There exists $\Lambda_\beta = \Lambda_\beta (M, g, k, \beta) \geq 1$ such that for every $\Lambda \geq \Lambda_\beta$ we have that $H({}_\Lambda b (|x|)) - \tr (k) ({}_\Lambda b(|x|)) > 0$ and $H( - {}_\Lambda b (|x|)) - \tr (k) (- {}_\Lambda b(|x|)) < 0$ on $\{x \in M : |x| > \Lambda\}$. We let $b_\beta(x) := {}_{\Lambda_\beta} b (|x|)$ on $\{x \in M : |x| > \Lambda_\beta \}$ and $c_\beta:= c \Lambda_\beta^{\beta -1}$.

\begin{proof} The Euclidean mean curvature of $\graph ({}_\Lambda b (|x|), \{x \in \R^n : |x| > \Lambda \})$ computed as the tangential divergence of the downward pointing unit normal equals $- \Lambda^{\beta-1} (n-\beta) |x|^{-\beta}$, and the length squared of its second fundamental form is $\Lambda^{2(\beta-1)} |x|^{- 2 \beta} (\beta^2 - 2 \beta + n)$. Since $(M,g, k)$ is asymptotically flat, we have that $ g_{ij} - \delta_{ij} = O^{1, 1 - n/p}(|x|^{-q})$. It follows that the difference between the mean curvature of this graph with respect to $g + dt^2$ and its Euclidean mean curvature is bounded above by a multiple of $\Lambda^{\beta -1} |x|^{- q - \beta}$ that depends only on $(M, g)$. Using that $k_{ij} = O^{0, 1 - n/p} (|x|^{-q-1})$, it follows that $\tr_g(k)(\pm {}_\Lambda b(|x|)) = O(|x|^{- \beta})$. The assertions follow easily from these estimates.   
\end{proof}
\end{proposition}

{\bf Assumptions: } For the remainder of this section we assume that $(M, g, k)$ is an asymptotically flat initial data set. If $n=3$, we assume that there exists $\beta \in (2, 3)$ such that $\tr_g (k) = O(|x|^{- \beta})$. (Note that if $(M, g, k)$ has harmonic asymptotics, then $\tr_g (k) = O(|x|^{-n})$, so that any $\beta \in (2, 3)$ will work.) If $n > 3$, we let $\beta := 1 + q$. Note that $\beta \in (\frac{n}{2}, n-1) \subset (2, n)$ and  that $\tr_g (k) = O(|x|^{- \beta})$. \\

Let $C := 1 +n \sup_{x \in M} |k(x)|_g$. Given $\tau > 0$, the constant functions $- \frac{C}{\tau}$ and $\frac{C}{\tau}$ are, respectively, sub and super solutions for (\ref{eqn:Jangequation}). For every $R \geq 1$ sufficiently large, there exists a unique solution $f_{\tau, R} \in \C^{3, \alpha}( \{x \in M : |x| \leq R\})$ of (\ref{eqn:Jangequation}) that vanishes on the boundary, cf. \cite[Lemma 2.2]{Eichmair:2009-Plateau}. The maximum principle applied as in \cite[Proposition 3]{Schoen-Yau:1981-pmt2} shows that $f_{\tau, R}(x)$ lies between $\pm \frac{C}{\tau}$ on $\{x \in M : |x| \leq R\}$ and between $\pm b_\beta$ on $\{x \in M : \Lambda_\beta < |x| \leq R\}$, where $b_\beta$ is as in Proposition \ref{prop:barriers}. It follows that the length of the gradient of $f_{\tau, R}$ is a priori bounded on compact subsets of $M$, independently of $R$ (but depending on $\tau$), cf. \cite[Lemma 2.1]{Eichmair:2009-Plateau} (a mistake in the proof is corrected in \cite[Appendix A]{potpourri}). Standard compactness results for solutions of the prescribed mean curvature equation, cf. e.g. \cite[Chapter 16]{Gilbarg-Trudinger:1998}, show that we can take a subsequential limit as $R \to \infty$ to obtain a solution $f_\tau \in \C^{3, \alpha}_{loc} (M)$ of (\ref{eqn:Jangequation}) with the following properties: 

\begin{proposition} \label{prop:existenceftau} For every $\tau >0$, there exists a solution $f_\tau \in \C^{3, \alpha}_{loc}(M)$ of (\ref{eqn:Jangequation}) such that $|f_\tau| \leq \frac{C}{\tau}$ on all of $M$. Moreover, $|f_\tau(x)| \leq c_\beta |x|^{2-\beta}$ for all $|x| \geq \Lambda_\beta$, where $c_\beta>0$ and $\Lambda_\beta \geq1$ are as in Proposition \ref{prop:barriers}. 
\end{proposition}

The almost minimizing property of graphs of bounded mean curvature is most conveniently expressed in terms of the ``$\lambda$-minimizing" currents of F. Duzaar and K. Steffen \cite{Duzaar-Steffen:1993}. We refer the reader to Appendix A in \cite{Eichmair:2009-Plateau} for a summary of  standard properties of $\lambda$-minimizing boundaries that we use here, and for further references to the geometric measure theory literature. 

\begin{proposition} [Cf. \protect {\cite[Example A.1] {Eichmair:2009-Plateau}}] \label{prop:geometricpropertiesftau} For $\tau>0$ let $f_\tau$ be as in Proposition \ref{prop:existenceftau}. Then $\graph(f_\tau) \subset M \times \R$ is a $2C$-minimizing boundary in $M \times \R$.
\end{proposition}

The results on the subsequential geometric limits $\bar \Sigma$ of $\graph(f_\tau, M) \subset M \times \R$ as $\tau \searrow 0$ stated in the following proposition correspond to those of \cite[Proposition 4]{Schoen-Yau:1981-pmt2} in the case $n=3$. For a derivation that is valid for all dimensions $3 \leq n < 8$ and which relies on the almost minimizing property of $\graph(f_\tau, M) \subset M \times \R$, we refer to \cite[p. 568--569, Remark 4.1 and Lemma 2.3]{Eichmair:2009-Plateau}. 

\begin{proposition} \label{prop:limitanalysis} Let $3 \leq n < 8$. There exists a properly embedded boundaryless $\C^{3, \alpha}_{loc}$ hypersurface $\bar \Sigma \subset M \times \R$ with the following properties: 
\begin{enumerate} [(a)]
\item $\bar \Sigma$ is the boundary of an open set $\Omega$. We have that $H_{\bar \Sigma} - \tr_{\bar \Sigma}(k) = 0$, where the mean curvature scalar $H_{\bar \Sigma}$ of $\bar \Sigma$ is computed as the tangential divergence of the unit normal pointing out of $\Omega$. Moreover, $\bar \Sigma = \partial \Omega$ is a $2C$-minimizing boundary. 
\item $\bar \Sigma$ has finitely many connected components. Each component $\Sigma$ of $\bar \Sigma$ is either cylindrical of the form $\Sigma_0 \times \R$, where $\Sigma_0$ is a closed properly embedded $\C^{3, \alpha}$ hypersurface in $M$, or it is the vertical graph of a $\C^{3, \alpha}_{loc}$ function $f_\Sigma$ whose domain $U_\Sigma$ is an open subset of $M$. The function $f_\Sigma$ is a solution of the Jang equation $H(f_\Sigma) - \tr(f_\Sigma) = 0$ on its domain $U_\Sigma$. 
\item  The boundary of the domain $U_\Sigma$ of every graphical component $\Sigma = \graph(f_\Sigma, U_\Sigma)$ of $\bar \Sigma$ is a closed properly embedded $\C^{3, \alpha}$ hypersurface in $M$. In fact, $\partial U_\Sigma$ consists of two disjoint unions $\Sigma_0^+$ and $\Sigma_0^-$ of components such that $f_\Sigma (x) \to \pm \infty$ uniformly as $x \to \Sigma_0^{\pm}$ in $U_\Sigma$. We have that $H_{\Sigma_0^\pm} \mp \tr_{\Sigma_0^\pm}(k) = 0$, where the mean curvature is computed as the tangential divergence of the unit normal pointing out of $U_\Sigma$. In fact, if $t \to \pm \infty$, then the hypersurfaces $\graph(f_\Sigma - t, U_\Sigma) \subset M \times \R$ converge locally uniformly in $\C^{3, \alpha}$ to the cylinder $\Sigma_0^\pm \times \R$. 
\item \label{exteriorsolution} $\bar \Sigma$ contains a graphical component $\Sigma = \graph(f_\Sigma, U_\Sigma)$ such that $\{x \in M : |x| > \Lambda_\beta\} \subset U_\Sigma$. We have that $|f_\Sigma(x)| \leq c_\beta |x|^{2 - \beta}$ for $|x| > \Lambda_\beta$ and that $f_\Sigma \in \C^{3, \alpha}_{loc} \cap W^{3, p}_{1-q} (\{x \in M : |x| > 2 \Lambda_\beta\})$. 
\end{enumerate} 
\begin{proof} We only comment on the asserted order of regularity and decay to assist the reader. 

The compactness and regularity theory for almost minimizing boundaries gives that subsequential $\C^{1, \alpha}_{loc}$ limits of the graphs exist. Such limits satisfy the marginally outer trapped condition weakly. Using standard arguments, one concludes that the convergence and limit are in $\C^{3, \alpha}_{loc}$. To derive the rate of decay of $f_\Sigma$ asserted in (\ref{exteriorsolution}), we first note that the interior gradient estimates for the Jang equation \cite[Lemma 2.1]{Eichmair:2009-Plateau} show that $|d f_\Sigma (x)|_g \to 0$ as $|x| \to \infty$. We can rescale the coordinate balls $B(x, |x|/2)$ to unit size, and with it $f_\Sigma, k_{ij}, g_{ij}$, and the equation. The rescaled metric converges to the Euclidean metric in $W^{2, p}$ (and hence in $\C^{1, 1 - n/p}$ by Sobolev embedding) and the rescaled second fundamental form converges to zero in $W^{1, p}$ (and hence in $\C^{0, 1 - n/p}$). We can use \cite[Theorem 13.1]{Gilbarg-Trudinger:1998} to obtain uniform H\"older estimates for the gradient of the rescaled function $f$. (Note that the equation is independent of ``$z$" in the notation of \cite{Gilbarg-Trudinger:1998}, so that ``$K$" is bounded independently of $|x|$.) The $W^{2, p}_{1-q}$ estimate follows from writing the equation in non-divergence form and applying $L^p$ theory (as in \cite[Theorem 9.11]{Gilbarg-Trudinger:1998}). This argument can be repeated upon differentiating the equation to give the asserted $W^{3, p}_{1-q}$ estimate.  

We note that when $n=3$, decay rates for $f_\tau$ are derived in \cite[p. 250]{Schoen-Yau:1981-pmt2} from their parametric estimate. The capillarity term makes their argument slightly more complicated. 
\end{proof}
\end{proposition}   

\begin{corollary} \label{cor:decayf} Let $f_\Sigma : U_\Sigma \to \R$ be the solution of the Jang equation on the exterior domain $U_\Sigma$ of (\ref{exteriorsolution}) in Proposition \ref{prop:limitanalysis}. The metric $\bar g := g + df_\Sigma \otimes df_\Sigma$ on $U_\Sigma \subset M$ is complete and in $\C^{2, \alpha}_{loc}$. In the asymptotically flat coordinates $(x_1, \ldots, x_n)$ on $\{x \in M : |x| > 2 \Lambda_\beta\}$ we have that $\bar g_{ij} - g_{ij} \in W^{2, p}_{-2 q} (\{ x \in M : |x| > 2 \Lambda_\beta\})$ and that $R_g - R_{\bar g} \in W^{0, p}_{-2q-2}(\{ x \in M : |x| > 2 \Lambda_\beta\})$. 
\end{corollary}

The constant $c_n$ that appears in the statements of the following two propositions is defined as $c_n := \frac{n-2}{4(n-1)}$. We recall that $k$ is extended to $M \times \R$ trivially in the vertical direction; $h$ is the second fundamental form tensor of $\Sigma$, oriented such that its trace is the mean curvature.  

\begin{proposition} [Cf. \protect{\cite[p. 254--256]{Schoen-Yau:1981-pmt2}}] \label{prop:SchoenYauidentity} Let $f_\Sigma: U_\Sigma \to \R$ be as in (\ref{exteriorsolution}) of Proposition \ref{prop:limitanalysis}. Let $\phi \in \C^1(\Sigma)$ be such that $(\supp \phi) \setminus \graph (f_\Sigma, \{x \in M : |x| >  2 \Lambda_{\beta}\})$ is compact in $\Sigma$. Then 
\begin{eqnarray} \label{eqn:SchoenYauidentityI}
c_n  \int_\Sigma 2 (\mu - |J|_g) \phi^2 + |h-k|_{\bar g}^2 \phi^2 d \L^n_{\bar g} + \frac{n}{2(n-1)} \int_\Sigma |d \phi|^2_{\bar g} d \L^n_{\bar g} \\ 
\leq \int_\Sigma |d \phi|^2_{\bar g} + c_n R_{\bar g} \phi^2 d \L^n_{\bar g}. \nonumber
\end{eqnarray}
Let $\Sigma_0^1, \ldots, \Sigma_0^l$ denote the components of $\partial U_\Sigma$. If the strict dominant energy condition $\mu > |J|_g$ holds near $\partial U_\Sigma$, then each $(\Sigma_0^i, g|_{\Sigma_0^i})$ has positive Yamabe type. In fact, the spectrum of the operator $- \Delta_{g|_{\Sigma_0^i}} \phi_i+ c_{n} R_{g|_{\Sigma_0^i} }$ is positive on each component. 

\begin{proof} We adapt the arguments in \cite{Schoen-Yau:1981-pmt2} to general dimensions. Consider the $1$-form $\omega =  d \left( \log \sqrt{1 + |d f_\Sigma|^2_g} \right) - k \lfloor \frac{\nabla_g f_\Sigma}{\sqrt{1 + |d  f_\Sigma|^2_g}}$ on $U_\Sigma \times \R$. Let $X \in \Gamma (T \Sigma)$ be the tangential vector field that is $\bar g$-dual to $\omega|_\Sigma$. Since $f_\Sigma = O^{2, \alpha}(|x|^{1-q})$, we have that $X = O(|x|^{-2q-1}) = o(|x|^{1-n})$. The crucial ``Schoen-Yau identity" \cite [(2.25)] {Schoen-Yau:1981-pmt2} reads \begin{equation} \label{eqn:SYidentity} \mu - J \left( \frac{\nabla_g f_\Sigma}{\sqrt{1 + |d f_\Sigma|^2_g}}\right) = \frac{1}{2} R_{\bar g} - \frac{1}{2} |h - k|_{\bar g}^2 -  |X|_{\bar g}^2 +  \div_{\bar g} X.\end{equation} (We refer the reader to \cite[Section 3.6]{Andersson-Eichmair-Metzger:2010} for an alternative, systematic derivation of this identity.)
Here, we evaluate at points $x \in U_\Sigma$ on the left and at points $(x, f_\Sigma(x)) \in \Sigma$ on the right. The left-hand side is bounded below by $\mu  - |J|_g$. Multiplying the pointwise equality (\ref{eqn:SYidentity}) by $\phi^2$, integrating over $\Sigma$, integrating by parts, and using the pointwise estimate $- |X|^2_{\bar g} \phi^2 - X(\phi^2) \leq |d \phi|^2_{\bar g}$,  one obtains that
\begin{eqnarray} \label{eqn:aux1}
\int_\Sigma 2 (\mu - |J|_g) \phi^2 +  |h - k|_{\bar g}^2  \phi^2 d \L^n_{\bar g} \leq \int_{\Sigma} R_{\bar g} \phi^2 + 2 |d  \phi|_{\bar g}^2 d \L^n_{\bar g}. 
\end{eqnarray} 
Recall that $\Sigma = \graph(f_\Sigma, U_\Sigma)$ has ends that are $\C^{3, \alpha}$-asymptotic to $(\partial U_\Sigma \times \R, g|_{\partial U_\Sigma} + dt^2)$. Using the assumption that $2 (\mu - |J|_g) \geq \delta > 0$ near $\partial U_\Sigma$, we see that (\ref{eqn:aux1}) implies that for every $\phi \in \C_c^1(\partial U_\Sigma \times \R)$,
\begin{eqnarray} \label{eqn:aux2}
\delta \int_{\partial U_\Sigma \times \R}  \phi^2  d \L^n_{g|_{\partial U_\Sigma} + dt^2} \leq \int_{\partial U_\Sigma \times \R}  R_{g|_{\partial U_\Sigma} + dt^2} \phi^2 + 2 |d \phi|^2_{g|_{\partial U_\Sigma} + dt^2} d \L^n_{g|_{\partial U_\Sigma} + dt^2}. 
\end{eqnarray}
The same ``separation of variables" argument as in \cite[p.254--255]{Schoen-Yau:1981-pmt2} shows that the estimate (\ref{eqn:aux2}) implies
\begin{eqnarray} \label{eqn:aux3}
\delta \int_{\partial U_\Sigma} \zeta^2  d \L^{n-1}_{g|_{\partial U_\Sigma}} \leq \int_{\partial U_\Sigma }  R_{g|_{\partial U_\Sigma}} \zeta^2 + 2 |d \zeta|^2_{g|_{\partial U_\Sigma}} d \L^{n-1}_{g|_{\partial U_\Sigma}} \text{ for every } \zeta \in \C^1(\partial U_\Sigma). 
\end{eqnarray}
Using a test function $\zeta$ that is one on one component of $\partial U_\Sigma$ and vanishes on all the other components when $n=3$, and using that $2 c_{n-1} \leq 1$ in the case $n > 3$, we conclude from (\ref{eqn:aux3}) that each component of $(\partial U_\Sigma, g|_{\partial U_\Sigma})$ has positive Yamabe type. Cf. with the argument in \cite{Galloway-Schoen:2006}. 
\end{proof}
\end{proposition}

The proof of the perturbation of the asymptotically cylindrical ends of $(\Sigma, \bar g)$ to exact cylindrical ends in Proposition \ref{prop:normalizedSigma}, and the method of conformal ``darning" of these exact cylindrical ends leading to inequality (\ref{eqn:variationalinequalityhat}) below, generalize the discussion in \cite[p. 256--257]{Schoen-Yau:1981-pmt2} to the case $3 \leq n <8$. \\

\begin{proposition} [Cf. \protect{\cite[p. 256--257]{Schoen-Yau:1981-pmt2}}]\label{prop:normalizedSigma}  Let $f_\Sigma: U_\Sigma \to \R$ be as in (\ref{exteriorsolution}) of Proposition \ref{prop:limitanalysis}. Assume that $U_\Sigma \neq M$,  that the dominant energy condition $\mu \geq |J|_g$ holds on $U_\Sigma$, and that this inequality is strict near $\partial U_\Sigma$. For every sufficiently large number $t_0 >1$ that is a regular value for both $f_\Sigma$ and $- f_{\Sigma}$ there exists a complete Riemannian metric $\tilde g$ on $\Sigma = \graph(f_\Sigma, U_\Sigma) \subset M \times \R$ with the following properties: 
\begin{enumerate}
\item There is a compact set $K \subset \Sigma$ such that the complement of $K$ in $\Sigma$ has finitely many components $N, C_1, \ldots, C_l$. We have that $N = \Sigma \cap (M \times (-t_0, t_0)) \supset (\graph(f_\Sigma, \{x \in M : |x| > 2 \Lambda_{\beta}\}), \bar g)$ and $\tilde g|_N = \bar g|_N$. Each $(C_i, \tilde g)$ is isometric to a half-cylinder $(\Sigma_{0}^i \times (0, \infty), g|_{\Sigma_0^i} + dt_i^2)$, where $\Sigma_{0}^1, \ldots, \Sigma_0^l$ are the components of $\partial U_\Sigma$. The metric $\tilde g$ is uniformly equivalent to the metric $\bar g$ on all of $\Sigma$.
\item For every $\phi \in \C^1 (\Sigma)$ such that $(\supp \phi ) \cap (C_1 \cup \ldots \cup C_l)$ is compact we have that
\begin{eqnarray} \label{eqn:variationalinequality}
\frac{1}{2} \int_{\Sigma}  |d \phi|^2_{\tilde g} d \L^n_{\tilde g} + c_n \int_{N} |h - k|^2_{\bar g} \phi^2 d \L^n_{\bar g} \leq \int_{\Sigma} |d  \phi |^2_{\tilde g} + c_n R_{\tilde g} \phi^2 d \L^n_{\tilde g}. 
\end{eqnarray}
\end{enumerate}
We will refer to $N$ as the asymptotically flat end of $\Sigma$, and to $C_1, \ldots, C_l$ as its cylindrical ends. 
\end{proposition}
If $U_\Sigma = M$, we take $\tilde g := \bar g$. \\

Let $(C_i, \tilde g) = (\Sigma_{0}^i \times (0, \infty), g|_{\Sigma_0^i} + dt_i^2)$ be one of the exact cylindrical ends of $(\Sigma, \tilde g)$. Let $0 < \phi_i \in \C^{2, \alpha} (\Sigma_0^i)$ be the first eigenfunction of the operator  $- \Delta_{g|_{\Sigma_0^i}} + c_{n} R_{g|_{\Sigma_0^i} }$, so that  $- \Delta_{g|_{\Sigma_0^i}} \phi_i+ c_{n} R_{g|_{\Sigma_0^i} } \phi_i = \lambda_i \phi_i$ for some $\lambda_i >0$. Let $\Psi_{C_i} := e^{- \sqrt{\lambda_i} t_i} \phi_i$. The scalar curvature of the metric $\Psi_{C_i}^\frac{4}{n-2} (g|_{\Sigma_0^i} + dt_i^2)$ on $\Sigma_0^i \times (0, \infty)$ vanishes. Let $s_i := \frac{n-2}{2 \sqrt{\lambda_i}} e^{- \frac{2 \sqrt{\lambda_i} t_i}{n-2}}$. Note that $(C_i, \Psi_{C_i}^{\frac{4}{n-2}} (g|_{\Sigma_0^i} + dt_i^2))$ is isometric to $((0, \frac{n-2}{2 \sqrt{\lambda_i}}) \times \Sigma_0^i,  \phi_i^{\frac{4}{n-2}} (\frac{4 \lambda_i s_i^2}{(n-2)^2}g|_{\Sigma_0^i} + ds_i^2))$ hence in particular uniformly equivalent\footnote{Two metrics $g_1, g_2$ on a manifold $M$ are uniformly equivalent if there exists a positive constant $c\geq1$ such that $c^{-1}g_2 \leq g_1 \leq c g_2$.} to the metric cone
$((0, \frac{n-2}{2 \sqrt{\lambda_i}}) \times \Sigma_0^i, s_i^2 g|_{\Sigma_0^i} + ds_i^2)$. Fix a function $0 < \Psi \in\C^{2, \alpha}_{loc}(\Sigma)$ that is equal to one on $K \cup N$ and equal to $\Psi_{C_i}$ on the part of $C_i$ where $t_i \geq 1$. Let $\tilde{g}_\Psi := \Psi^{\frac{4}{n-2}} \tilde g$. The scalar curvature $R_{\tilde{g}_\Psi}$ of $\tilde{g}_\Psi$ vanishes on each cylindrical end of $\Sigma$. Replace $\phi$ by $\phi \Psi$ in  (\ref{eqn:variationalinequality}). Using that $- \Delta_{\tilde g} \Psi + c_n R_{\tilde g} \Psi = c_n R_{\tilde g_{\Psi}} \Psi^{\frac{n+2}{n-2}}$, that $\Psi$ is one on $N$, that $\phi$ is compactly supported in $\Sigma \setminus N$, and an integration by parts, we obtain that 
\begin{eqnarray} \label{eqn:variationalinequalityhat}
\int_{\Sigma} \frac{1}{2} \Psi^{-2} |d (\Psi \phi)|_{\tilde{g}_\Psi}^2 d \L^n_{\tilde{g}_\Psi} + c_n \int_{N} |h - k|^2_{\bar g} \phi^2 d \L^n_{\bar g}\leq \int_{\Sigma}  |d \phi |^2_{\tilde{g}_\Psi} + c_n R_{\tilde{g}_\Psi} \phi^2  d \L^n_{\tilde{g}_\Psi}
\end{eqnarray} 
for every  $\phi \in \C^1(\Sigma)$ such that $(\supp \phi ) \cap (C_1 \cup \ldots \cup C_l)$ is compact.\footnote{Note that $\Sigma$ is not complete with respect to the $\tilde{g}_\Psi$ metric unless $U_\Sigma = M$.  We will always refer to the topology on $\Sigma$ that comes from the complete metric $\bar g$.} \\

We introduce a new ``distance" function $0 < s \in \C^{3, \alpha}_{loc}(\Sigma)$ such that $s(x) = |x|$ when $x \in N$ and such that $s(x) = s_i(x)$ when $x \in C_i$ and $s_i (x) \leq \frac{n-2}{2 \sqrt{\lambda_i}} e^{- \frac{2 \sqrt{\lambda_i}}{n-2}}$. (When $U_\Sigma = M$ we take $s(x) = |x|$.) \\

The metric completion of $(\Sigma, \tilde{g}_\Psi)$ as a metric space is obtained by adding a point at infinity for each cylindrical end in $(\Sigma, \tilde g)$. Note that $s(x)$ would extend continuously (by zero) to these points in the completion. It is helpful to think of these points as virtual singular points of $(\Sigma, \tilde{g}_\Psi)$. The $\tilde{g}_\Psi$--harmonic capacity of these singular points vanishes. More precisely, for $\epsilon \in (0, 1)$ small consider the functions $\chi_\epsilon \in \C^{3, \alpha}_{loc}(\Sigma)$ defined by 
\begin{displaymath}
\chi_\epsilon (x) := \left\{ \begin{array}{ll} 0 & 0 < s(x) \leq \epsilon \\ - 1 + \epsilon ^{-1} s(x) & \epsilon \leq s(x) \leq 2 \epsilon \\ 1 & 2 \epsilon \leq s(x). \end{array} \right.
\end{displaymath}
Note that $\chi_\epsilon \equiv 1$ on $K \cup N$, that $\chi_\epsilon$ has compact support in $\Sigma \setminus (K \cup N) = C_1 \cup \ldots \cup C_l$, that $\chi_\epsilon \to 1$ locally uniformly on all of $\Sigma$ as $\epsilon \searrow 0$, and that $\int_{\Sigma} |d \chi_\epsilon|^2_{\tilde{g}_\Psi} d \L^n_{\tilde{g}_\Psi} = O(\epsilon^{n-2})$. 


\section{The conformal structure of $(\Sigma, \tilde{g}_\Psi)$} \label{sec:conformalstructure}

When $(M, g, k)$ has harmonic asymptotics, the decay of $f_\Sigma$ and the properties of the metric $\bar g$ on $\Sigma$ in Corollary \ref{cor:decayf} can be improved:

\begin{proposition} \label{prop:betterasymptotics} Let $(M, g, k)$ be asymptotically flat with harmonic asymptotics of type $(2, \alpha)$. For every $\eta >0$ sufficiently small, we can arrange for the decay of $f_\Sigma$ in Proposition \ref{prop:limitanalysis} (\ref{exteriorsolution}) to be $f_\Sigma = O^{3, \alpha} (|x|^{2 - n + \eta})$. The assertions of Corollary \ref{cor:decayf} are improved to $\bar g_{ij} = g_{ij} + O^{2, \alpha}(|x|^{2-2n+2\eta})$ and $R_{\bar g} = R_{g} + O^{0, \alpha} (|x|^{-2 n + 2 \eta})$. 
\end{proposition}

{\bf Assumptions: } We assume throughout this section that $(M, g, k)$ has harmonic asymptotics of type $(2, \alpha)$ such that $\mu > |J|_g$ and such that $(\mu, J) = O^{0, \alpha} (|x|^{-n-q_0'})$ for some $q_0' > 1$. In particular, we have that $\bar g_{ij} = g_{ij} + O^{2, \alpha}(|x|^{2 - 2 n + \delta})$ and that $R_{\bar g} = O^{0, \alpha}(|x|^{-n-1-\delta})$ for all $\delta >0$ sufficiently small. We will work with the (typically) incomplete Riemannian manifold $(\Sigma, \tilde{g}_\Psi)$ associated with $(M, g, k)$ that has been constructed in Section \ref{sec:geometricsolutions}, and the modified distance function $s \in \C^{2, \alpha}(\Sigma)$.  \\


The following proposition corresponds to Lemma 4 in \cite{Schoen-Yau:1981-pmt2}. The difference here is that $(\Sigma, \tilde{g}_\Psi)$ is in general not uniformly equivalent to a smooth manifold. Apart from necessary technical modifications, the proof is very similar to that of \cite{Schoen-Yau:1981-pmt2}. \\

\begin{proposition} [Cf. \protect{\cite[Lemma 4]{Schoen-Yau:1981-pmt2}}] \label{prop:conformalfactor} There exists $u \in \C^{2, \alpha}_{loc}(\Sigma)$ such that for some $c \geq 1$ one has that $c^{-1} \leq u(x) \leq c$, such that $- \Delta_{\tilde{g}_\Psi} u + c_n R_{\tilde{g}_\Psi} u = 0$ on $\Sigma$, and such that 
\begin{eqnarray*}
u (x) - 1 - a_1 |x|^{2-n} = O^{2, \alpha} (|x|^{1-n}) \text{ as } |x| \to \infty.  
\end{eqnarray*} 
for some constant $a_1 \leq 0$. 
\begin{proof} Fix $\sigma_0>0$ small so that the scalar curvature of $\tilde{g}_\Psi$ vanishes on $ \{ s(x) < 2 \sigma_0\}$, and such that for all $\sigma \in (0, 2 \sigma_0)$ both $\sigma$ and $\sigma^{-1}$ are regular values for $s(x)$. For $\sigma \in (0, \sigma_0)$, consider the Dirichlet problems 
\begin{displaymath}
\left\{ \begin{array}{lllll} - \Delta_{\tilde{g}_\Psi} v_\sigma + c_n R_{\tilde{g}_\Psi}  v_\sigma &=& - c_n R_{\tilde{g}_\Psi} &\text{ on }& \{\sigma < s(x) < \sigma^{-1}\} \\ v_\sigma &=& 0 &\text{ on }&  \{s(x) = \sigma\} \cup \{s(x) = \sigma^{-1}\}. \end{array} \right.
\end{displaymath}
In view of (\ref{eqn:variationalinequalityhat}), the corresponding homogeneous problems (with zero right-hand side) only admit the zero solution. Hence, by the Fredholm alternative, the above Dirichlet problems admit unique solutions $v_\sigma$. We extend each of these solutions $v_\sigma$ by $0$ to a compactly supported Lipschitz function on all of $\Sigma$. 

We have that 
\begin{eqnarray*} 
&& \left( \int_{\{s(x) \geq \sigma_0\}} |v_\sigma|^{\frac{2n}{n-2}} d \L_{\tilde{g}_\Psi}^n\right)^{\frac{n-2}{n}}  \leq C_1 \left( \int_{\{s(x) \geq \sigma_0\}} |\Psi v_\sigma|^{\frac{2n}{n-2}} d \L_{\tilde{g}_\Psi}^n \right)^{\frac{n-2}{n}} \\ &\leq& C_2 \int_{\{s(x) \geq \sigma_0\}} \Psi^{-2} |d (\Psi v_\sigma)|^2_{\tilde{g}_\Psi} d \L_{\tilde{g}_\Psi}^n \leq C_2 \int_{\{\sigma \leq s(x) \leq \sigma^{-1}\}} \Psi^{-2} |d (\Psi v_\sigma)|_{\tilde{g}_\Psi}^2  d \L_{\tilde{g}_\Psi}^n \\ &\leq& 2 C_2 \int_{\{\sigma  \leq s(x) \leq \sigma^{-1}\}} |R_{\tilde{g}_\Psi}| |v_\sigma| d \L_{\tilde{g}_\Psi}^n \\ &\leq& 2 C_2 \left( \int_{\{s(x) \geq \sigma_0\}} |R_{\tilde{g}_\Psi}|^{\frac{2n}{n+2}} d \L_{\tilde{g}_\Psi}^n \right)^{\frac{n+2}{2n}}  \left( \int_{\{s(x) \geq \sigma_0\}} |v_\sigma|^{\frac{2n}{n-2}} d \L_{\tilde{g}_\Psi}^n\right) ^{\frac{n-2}{2n}}
\end{eqnarray*}
where $C_1, C_2>0$ are constants that do not depend on $\sigma \in (0, \sigma_0)$. The first inequality holds because $\Psi$ is bounded below on $\{s(x) \geq \sigma_0\}$, the second follows from the Sobolev inequality in the form of Lemma \ref{lem:SobolevInequality} on $(\{s(x) \geq \sigma_0\}, \tilde{g}_\Psi)$ and because $\Psi$ is bounded above, the third follows from inclusion, the forth from (\ref{eqn:variationalinequalityhat}) (multiply the equation that $v_\sigma$ satisfies by $v_\sigma$ and integrate by parts), and the last inequality from H\"older's inequality and the fact that $R_{\tilde{g}_\Psi}$ is supported in $\{s(x) \geq \sigma_0\}$. Note that $\int_{\{s(x) \geq \sigma_0\}} |R_{\tilde{g}_\Psi}|^{\frac{2n}{n+2}} d \L_{\tilde{g}_\Psi}^n < \infty$ because $R_{\tilde{g}_\Psi} = R_{\bar g} = R_{\tilde g} \in \C^{0, \alpha}_{-n-1}$ on $N$. We conclude that $\int_{\{s(x) \geq \sigma_0\}} |v_\sigma|^{\frac{2n}{n-2}} d \L_{\tilde{g}_\Psi}^n $ is bounded independently of $\sigma \in (0, \sigma_0)$. Standard elliptic theory shows that $v_\sigma$ is uniformly bounded in $\C^{2, \alpha}$ on $\{s(x) \geq  2 \sigma_0\}$. Since $v_\sigma$ is harmonic on $\{\sigma < s(x) < 2 \sigma_0\}$ and because harmonic functions achieve their maximum and their minimum values on the boundary, it follows that $v_\sigma$ is uniformly bounded above and below on all of $\{\sigma < s(x) < \sigma^{-1}\}$. Standard theory converts the $L^\infty$-estimate into $\C^{2, \alpha}_{loc}$ bounds on all of $\{\sigma < s(x) < \sigma^{-1}\}$. \\

Let $u_\sigma := v_\sigma + 1$, and let $u \in \C^{2, \alpha}_{loc}$ denote a subsequential limit of $\{u_\sigma\}_{\sigma \in (0, \sigma_0)}$ as $\sigma \searrow 0$. We first claim that $u_\sigma > 0$ on $\{\sigma < s(x) < \sigma^{-1}\}$. Clearly this is true near the boundary of $\{\sigma < s(x) < \sigma^{-1}\}$. Let $\epsilon > 0$ small be a regular value of $- u_\sigma$. Note that $\min\{u_\sigma + \epsilon, 0\}$ is a Lipschitz function with support in $\{\sigma < s(x) < \sigma^{-1}\}$. Using it as a test function in (\ref{eqn:variationalinequalityhat}) we obtain that 
\begin{eqnarray*} 
&& \frac{1}{2} \int_{\{u_\sigma < - \epsilon\}}    \Psi^{-2} |d (\Psi (u_\sigma + \epsilon))|^2_{\tilde{g}_\Psi} d \L^n_{\tilde{g}_\Psi} \\ 
&\leq& \int_{\{u_\sigma < - \epsilon\}} |d  (u_\sigma + \epsilon)|^2_{\tilde{g}_\Psi} + c_n R_{\tilde{g}_\Psi} (u_\sigma + \epsilon)^2 d \L^n_{\tilde{g}_\Psi} \\
&=& \int_{\{u_\sigma < - \epsilon\}} (u_\sigma + \epsilon) (- \Delta_{\tilde{g}_\Psi} (u_\sigma + \epsilon)  + c_n R_{\tilde{g}_\Psi} (u_\sigma + \epsilon)) d \L^n_{\tilde{g}_\Psi} \\ &=& \int_{\{u_\sigma < - \epsilon\}} c_n \epsilon (u_{\sigma} + \epsilon) R_{\tilde{g}_\Psi} d \L^n_{\tilde{g}_\Psi}. 
\end{eqnarray*}
Letting $\epsilon \searrow 0$, we see that $\Psi u_\sigma$ is constant on $\{u_\sigma<  0\}$. It follows that $\{u_\sigma < 0\}$ is empty and hence that $u_\sigma \geq 0$ on $\{\sigma < s(x) < \sigma^{-1}\}$. By Harnack theory, $u_\sigma >0$ on $\{\sigma < s(x) < \sigma^{-1}\}$. \\

Recall that $R_{\tilde{g}_\Psi} = R_{\bar g} = O^{0, \alpha} (|x|^{-n-1-\delta})$ for $|x|$ large. In conjunction with the bound for $||v_{\sigma}||_{L^{2n/(n-2)} (\{s(x) \geq \sigma_0\})}$ obtained above, standard elliptic theory shows that $u(x) \to 1$ as $|x| \to \infty$. Harnack theory gives that $u>0$ on all of $\Sigma$. Using that $u_\sigma$ is harmonic on $\{\sigma < s(x) < 2 \sigma_0\}$, we conclude from the maximum principle that the $u_\sigma$ are bounded uniformly above and below on $\{\sigma < s(x) < \sigma^{-1}\}$ by positive constants  that are independent of $\sigma \in (0, \sigma_0)$.  Thus $u$ is bounded above and below by positive constants on all of $\Sigma$.  The existence of $a_1 \in \R$ such that $u - 1 - a_1 |x|^{2-n} = O^{2, \alpha} (|x|^{1-n})$ follows from asymptotic analysis as in e.g. \cite{Bartnik:1986}.\footnote{ \label{footnote:aux} When $n=3$ we use that $\tilde g_\psi = \varphi^{4} (\delta_{ij} + O^{2, \alpha} (|x|^{-3}))$ with $\varphi (x) = 1 + a |x|^{-1} + O^{2, \alpha} (|x|^{-2})$ as well as the conformal invariance of the equation $- 8 \Delta_{\tilde g_\Psi} u + R_{\tilde g_\Psi} u = 0$ to obtain this expansion.} \\

Our next goal is to show that $\int_\Sigma |d u|^2_{\tilde{g}_\Psi} d \L^n_{\tilde{g}_\Psi} < \infty$. Since $- \Delta_{\tilde{g}_\Psi} u + c_n R_{\tilde{g}_\Psi} u = 0$ and $R_{\tilde{g}_\Psi} \equiv 0$ on $\{s(x) \leq \sigma_0\}$, we have that 
$$\int_{\{s(x) < \sigma_0\}} d u ( \nabla_{\tilde{g}_\Psi} \xi)  d \L^{n}_{\tilde{g}_\Psi}  = \int_{\{s(x) = \sigma_0\}} \xi \nu_{\tilde{g}_\Psi}(u) d \H^{n-1}_{\tilde{g}_\Psi} \text{ for every } \xi \in \C^1_c( \{s(x) \leq \sigma_0\}).$$
Applying this identity with the test functions $\xi = u \chi_\epsilon^2$, letting $\epsilon \searrow 0$, and using our $L^\infty$-bound for $u$, we obtain that $\int_{\{s(x) \leq \sigma_0\}} |d  u|_{\tilde{g}_\Psi}^2 d \L^n_{\tilde{g}_\Psi} < \infty$. Using the test functions $\xi = u \chi_\epsilon$ and letting $\epsilon \searrow 0$, we conclude that for every $\sigma \in (0, \sigma_0)$ 
\begin{eqnarray} \label{eqn:Greensidentity} \int_{\{s(x) < \sigma\}} |d u|_{\tilde{g}_\Psi}^2  d \L^{n}_{\tilde{g}_\Psi}  = \int_{\{s(x) = \sigma\}} u \nu_{\tilde{g}_\Psi}(u) d \H^{n-1}_{\tilde{g}_\Psi}.\end{eqnarray} Also, since $d u = O^{1, \alpha} (|x|^{1-n})$, we obtain that $\int_{\Sigma} |d u |^2_{\tilde{g}_\Psi} d \L^n_{\tilde{g}_\Psi} < \infty$. \\

Finally, to get that $a_1 \leq 0$, we use (\ref{eqn:Greensidentity}), (\ref{eqn:variationalinequalityhat}), the finiteness of $\int_\Sigma |d u|_{\tilde{g}_\Psi}^2 d \L^n_{\tilde{g}_\Psi}$, the integrability of $|h - k|^2_{\bar g}$ on $N$, an integration by parts, and that $- \Delta_{\tilde g_\Psi} u + c_n R_{\tilde g_\Psi} u = 0$ to obtain that 
\begin{eqnarray} \label{eqn:a1}
0 &\leq& \int_N \frac{1}{2} |du|_{\bar g}^2 + c_n |h - k|^2_{\bar g} u^2 d \L^n_{\bar g} \\  &\leq& \nonumber \liminf_{\epsilon \searrow 0} \int _\Sigma |d (\chi_\epsilon u)|_{\tilde{g}_\Psi}^2 + c_n R_{\tilde{g}_\Psi} (\chi_\epsilon u)^2 d \L^n_{\tilde{g}_\Psi} \\ 
&=& \int_\Sigma |d u|_{\tilde{g}_\Psi}^2 + c_n R_{\tilde{g}_\Psi} u^2 d \L^n_{\tilde{g}_\Psi} \nonumber 
= \lim_{\sigma \searrow 0}  \int_{\{\sigma < s(x) < \sigma^{-1}\}} |d u|_{\tilde{g}_\Psi}^2 + c_n R_{\tilde{g}_\Psi} u^2 d \L^n_{\tilde{g}_\Psi} \nonumber \\
&=& \lim_{\sigma \searrow 0} \int_{\{s(x) = \sigma^{-1}\}} u \nu_{\bar g} (u) d \H^{n-1}_{\bar g} - \lim_{\sigma \searrow 0} \int_{\{ s(x) \leq \sigma \}} |d u|_{\tilde g_\Psi}^2 d \L^n_{\tilde g_\Psi}  \nonumber \\
&=& (2-n) a_1 |\mathbb{S}^{n-1}| \nonumber
\end{eqnarray}
(Recall that $\bar g = \tilde{g}_\Psi$ on the asymptotically flat end $N$ of $\Sigma$.) 
\end{proof}
\end{proposition} 
Following \cite{Schoen-Yau:1981-pmt2}, we define a new metric $\tilde{g}_{u\Psi} := u^{\frac{4}{n-2}} \tilde{g}_\Psi = (u \Psi)^{\frac{4}{n-2}} \tilde g$ on $\Sigma$. Note that $R_{\tilde{g}_{u\Psi}} = 0$ since $- \Delta_{\tilde g_{\Psi}} u + c_n R_{\tilde g_\Psi } u = 0$. The following proposition, whose proof is similar to that of Proposition \ref{prop:conformalfactor}, undoes the conformal darning of the cylindrical ends in $(\Sigma, \tilde g)$ that was effected by the conformal factor $\Psi$. It is our substitute for the explicit Green's function in \cite[p. 259]{Schoen-Yau:1981-pmt2} with poles at the singular points in dimension $n=3$.

\begin{proposition} \label{prop:blowupfunction}
There exists $0 < w \in \C^{2, \alpha}_{loc} (\Sigma)$ such that $\Delta_{\tilde{g}_{u\Psi}} w \leq 0$ with strict inequality when $|x|$ is large, such that $w(x) - a_2 |x|^{2-n} = O^{2, \alpha} (|x|^{1-n})$ as $|x| \to \infty$ for some constant  $a_2 \in \R$, and such that for some $c \geq 1$, one has that $\frac{1}{c} \frac{1}{u s^{n-2}} \leq w(x) \leq  \frac{c}{ u s^{n-2}}$ as $s(x) \to 0$. 
\begin{proof} 
Let $\sigma_0>0$ be as in the proof of Proposition \ref{prop:conformalfactor}. Note that $\Delta_{\tilde{g}_{u\Psi}} \frac{1}{u s^{n-2}} = 0$ on $\{s(x) < 2 \sigma_0\}$. Fix a non-negative function $w_0 \in \C^{2, \alpha}_{loc}(\Sigma)$ that agrees with $\frac{1}{u s^{n-2}}$ on  $\{s(x) < 2 \sigma_0\}$ and such that $\supp (w_0) \cap \{s(x) > \sigma_0\}$ is compact.  Fix a non-negative function  $q \in \C^{2, \alpha}_{loc} (\Sigma)$ with $\supp (q) \cap \{ s(x) < 2 \sigma_0\} = \emptyset$ and such that $q(x) = |x|^{-2n}$ when $|x|$ is large. Given $\sigma \in (0, \sigma_0)$, let $w_\sigma$ be the unique solution of
\begin{displaymath}
\left\{ \begin{array}{lllll} \Delta_{\tilde{g}_{u\Psi}} (w_0 + w_{\sigma}) &=& - q &\text{ on }& \{\sigma < s(x) < \sigma^{-1}\} \\ w_\sigma &=& 0 &\text{ on }& \{s(x) = \sigma\} \cup \{s(x) = \sigma^{-1}\}. \end{array} \right.
\end{displaymath} 
Note that $w_0 + w_{\sigma}$ is positive by the maximum principle. We extend $w_{\sigma}$ by $0$ to a Lipschitz function on all of $\Sigma$. We have that
\begin{eqnarray*}
&& C_1^{-1}  \left(  \int_{\{s(x) \geq \sigma_0\} }  |w_\sigma|^{\frac{2n}{n-2}} d \L^n_{\tilde{g}_{u\Psi}}\right)^{\frac{n-2}{n}} \\ &\leq& \int_{\{s(x) \geq \sigma_0\}} |d  w_{\sigma}|^2_{\tilde{g}_{u\Psi}}  d \L^n_{\tilde{g}_{u\Psi}} \leq \int_{\{\sigma \leq s(x) \leq \sigma^{-1}\}} |d w_\sigma|^2_{\tilde{g}_{u\Psi}} d \L^n_{\tilde{g}_{u\Psi}} \\ &=& \int_{\{\sigma \leq s(x) \leq \sigma^{-1}\}} w_\sigma (q + \Delta_{\tilde{g}_{u\Psi}} w_0) d \L^n_{\tilde{g}_{u\Psi}} \leq \int_{\{s(x) \geq \sigma_0\}} |w_\sigma| |q+ \Delta_{\tilde{g}_{u\Psi}} w_0| d \L^n_{\tilde{g}_{u\Psi}} \\ &\leq& \left( \int_{\{s(x) \geq \sigma_0\}} |w_\sigma|^{\frac{2n}{n-2}} d \L^n_{\tilde{g}_{u\Psi}}\right)^{\frac{n-2}{2n}} \left( \int_{\{s(x) \geq \sigma_0\}} |q + \Delta_{\tilde{g}_{u\Psi}} w_0|^{\frac{2n}{n+2}} d \L^n_{\tilde{g}_{u\Psi}}\right)^{\frac{n+2}{2n}}. 
\end{eqnarray*}
We have used the Sobolev inequality in the form of Lemma \ref{lem:SobolevInequality} on $(\{s(x) \geq \sigma_0\}, \tilde{g}_{u\Psi})$ in the first inequality. It follows that $\int_{\{s(x) \geq \sigma_0\}} |w_\sigma|^{\frac{2n}{n-2}} d \L^n_{\tilde{g}_{u\Psi}}$ is bounded independently of $\sigma \in (0, \sigma_0)$. From this and the equation that $w_\sigma$ satisfies we obtain  $\C^{2, \alpha}$ estimates for $w_\sigma$ on $\{s(x) \geq 2 \sigma_0\}$. Using that $w_\sigma$ is $\tilde{g}_{u\Psi}$-harmonic on $\{ \sigma < s(x) < 2 \sigma_0\}$ and that $w_\sigma$ vanishes on $\{s(x) = \sigma\}$, we obtain an $L^\infty$-bound for $w_\sigma$ that is independent of $\sigma \in (0, \sigma_0)$. Passing to a subsequential limit $\sigma_i \searrow 0$, we obtain a non-negative function $w:= w_0 + \lim_{i \to \infty} w_{\sigma_i} \in \C^{2, \alpha}_{loc}(\Sigma)$ such that $\Delta_{\tilde{g}_{u\Psi}} w = - q$. Using the Harnack principle, standard asymptotic analysis (see also footnote \ref{footnote:aux}), the $L^\infty$-bound for $w - w_0$, and that $u$ is bounded above and below by positive constants, we see that $w$ has all the asserted properties.  
\end{proof}
\end{proposition}


\section{Proof that $E \geq 0$} \label{sec:egeq0}

The argument in Proposition \ref{prop:main} below follows \cite[p. 259]{Schoen-Yau:1981-pmt2}. We supply additional details to explain why the Riemannian positive energy theorem can be applied for time-symmetric initial data sets with certain non-standard ends.
\begin{proposition} \label{prop:main} Assumptions as in the first part of Theorem \ref{thm:main}. Then $E \geq 0$. 
\begin{proof}

By Theorem \ref{thm:density}, every asymptotically flat initial data set satisfying the hypotheses of the first part of Theorem \ref{thm:main} can be approximated by initial data sets that satisfy the assumptions of Section \ref{sec:conformalstructure} and whose energies converge to the energy of $(M,g, k)$. We may and will assume that $(M, g, k)$ is such a special initial data set. \\

In view of (\ref{eqn:EPab}), we need to show that $a \geq 0$, where $a$ is as in (\ref{eqn:uY}). We let $(\Sigma, \tilde g)$, $\Psi$, $u$, and $w$ as in Proposition \ref{prop:normalizedSigma}, the discussion succeeding it,  Proposition \ref{prop:conformalfactor}, and Proposition \ref{prop:blowupfunction} respectively. \\

Fix $\epsilon >0$ small. Define a new metric $\tilde g^\epsilon := \left((1 + \epsilon w) u \Psi \right)^{{4}/{(n-2)}} \tilde g$ on $\Sigma$. By Proposition \ref{prop:conformalfactor}, the metric $\tilde{g}_{u\Psi} = ( u\Psi )^{4/(n-2)} \tilde g$ on $\Sigma$ is scalar flat. Proposition \ref{prop:blowupfunction} shows that the scalar curvature of $\tilde g^\epsilon$ is non-negative and positive when $|x|$ is large. In the asymptotically flat coordinate chart $\tilde g^\epsilon$ has the expansion
$$\tilde g^\epsilon_{ij}  =  \left( 1 +  \frac{(n-2)(a + a_1  + \epsilon a_2)}{4}  |x|^{2-n} \right)^{\frac{4}{n-2}} \delta_{ij} + O^{2, \alpha} (|x|^{1-n}) \text{ as } |x| \to \infty.$$ 
Recall that $\{ s(x) < \sigma_0\}$ has $l$ connected components, one for each of the cylindrical ends $C_1, \ldots, C_l$ of $(\Sigma, \tilde g)$. On each of these components, $\tilde g^\epsilon$ is uniformly equivalent to a metric of the form $\sigma^2_i  g|_{\Sigma_0^i} + d \sigma_i^2$, where $\sigma_i (x):= s(x)^{- (n-2)/4}$ on $C_i$. It follows that $\tilde g^\epsilon$ is a complete metric on $\Sigma$. As in \cite[p. 259]{Schoen-Yau:1981-pmt2}, we conclude from the minimal hypersurface proof of the Riemannian positive energy theorem \cite{Schoen-Yau:1979-pmt1, Schoen-Yau:1981-asymptotics, Schoen-Yau:1979-pat, Schoen:1989} that \begin{eqnarray} \label{ineq:ea_1} a  + a_1 + \epsilon a_2 \geq 0.\end{eqnarray} Since this holds for every $\epsilon >0$, and since $a_1 \leq 0$, we conclude that indeed $a \geq 0$. \\

Non-asymptotically flat ends in $(\Sigma, \tilde g^\epsilon)$ such as $(C_i, \tilde g^\epsilon)$ were not considered in the original statement and proof of the Riemannian positive energy theorem \cite{Schoen-Yau:1979-pmt1, Schoen:1989}, so we include two additional remarks that explain why we can use this result here. \\

The first step in the minimal hypersurface proof of the Riemannian positive energy theorem is a deformation from non-negative scalar curvature to everywhere positive scalar curvature. What is actually required later in the proof is that the complete connected area minimizing hypersurface that is constructed as a limit of least area surfaces spanning the co-dimension $2$ surfaces $\{(x_1, \ldots, x_{n-1}, h_\rho) : \sum_{i=1}^{n-1} x_i^2 = \rho^2\}$ where $h_\rho \in [- \Lambda, \Lambda]$ for a certain $\Lambda \geq 1$, as $\rho \to \infty$, passes through some region where the scalar curvature is strictly positive. Since this hypersurface is unbounded in the asymptotically flat end of $(\Sigma, \tilde g^\epsilon)$, and since $R_{ \tilde g^\epsilon} >0$ when $|x|$ is large, this will always be the case. \\

One must also ensure that the area minimizing hypersurface that is constructed intersects the non-asymptotically flat ends $C_1, \ldots, C_l$ of $(\Sigma, \tilde g^\epsilon)$ in a compact set. To see this, one argues that the hypersurface measure of the part of the minimizing hypersurface in $\{x \in C_i : \sigma_ i (x) \in (\sigma, \sigma+1)\}$ -- if non-empty -- is bounded below by a positive constant that is independent of $\sigma$. We cannot appeal to the usual monotonicity formula for stationary hypersurfaces directly, because we don't have sufficient control on the metric $\tilde g^\epsilon$ as $\sigma_i (x) \to \infty$. However, a variant of a classical argument of De Giorgi is available in the form of Lemma \ref{lem:monotonicity} to argue this, cf \cite[Proof of Lemma 3.3]{Bray-Lee:2009}. If the connected minimizing hypersurface reached too far into one of the ends $C_i$, then one could use this estimate to contradict the area minimizing property by capping off these ``fingers". \\
\end{proof}
\end{proposition}


\section{Rigidity: The case $E = 0$}

The following proposition is a simple consequence of the proof of Theorem \ref{thm:density} from \cite{spacetimePMT}:

\begin{proposition} \label{prop:densitytrace} Hypotheses as in Theorem \ref{thm:density}. Then $(g^j, k^j) \to (g, k)$ in $\C^{2, \alpha}_{loc} \times \C^{1, \alpha}_{loc}$. When $n=3$ and if there exists $\gamma > 2$ such that $\tr_g(k) =  O(|x|^{- \gamma})$, then, given $\gamma' \in (2, \min\{\gamma, q_0 + 2, 2 q +1\})$, the approximating initial data $(g^j, k^j)$ in the conclusion of Theorem \ref{thm:density} can be chosen such that $\tr_{g^j} (k^j) = O(|x|^{- \gamma'})$ uniformly in $j$.  
\begin{proof} 
We remark that the vector fields $Z$ in Lemma 21 and $Y$ in Lemma 23 of \cite{spacetimePMT} have expansions of the form $c_i |x|^{2-n} + O^{1, \alpha} (|x|^{2 - n - \eta})$, where the constants $c_i$ and the error term depend continuously on $(g, \pi) \in W^{2, p}_{-q} \times W^{1, p}_{-q-1}$ and $(\mu, J) \in \C^{0, \alpha}_{-n-q_0}$, and the parameters of the construction for any $\eta \in (0, \min\{2q+2-n, q_0\})$. Compare with equation (37) and the proof of Proposition 24 in \cite{spacetimePMT}.     
\end{proof}
\end{proposition}

The main idea of the proof below follows \cite[Section 6]{Schoen-Yau:1981-pmt2}. Our regularity and decay assumptions on the initial data set are weaker than those of \cite{Schoen-Yau:1981-pmt2}. Moreover, the potential cylindrical ends cause some additional complications in dimensions $3 < n < 8$. This is why we give a detailed proof here. The use of the Cheeger-Gromoll splitting theorem (to rule out cylindrical ends) and the Bishop-Gromov volume comparison theorem (to show that the Jang graph is isometric to Euclidean space) in the final step  to obtain rigidity in the Riemannian positive energy theorem is also different from the classical argument in e.g. \cite[Proposition 2]{Schoen:1984}. 

\begin{proposition}[Cf. \protect{\cite[Section 6]{Schoen-Yau:1981-pmt2}}]  \label{prop:rvanishes} Assume that the asymptotically flat initial data set $(M, g, k)$ satisfies all the conditions of Theorem \ref{thm:main} including the additional hypothesis that $\tr_g(k) = O(|x|^{-\gamma})$ for some $\gamma >2$ when $n=3$. If $E = 0$, then $(M, g, k)$ is Cauchy initial data for Minkowski space $(\R^n \times \R, dx_1^2 + \ldots + dx_n^2 - dx_{n+1}^2)$. 
\begin{proof}
Let $(g^j, k^j)$ be a sequence of data on $M$ of harmonic type $(2, \alpha)$ that approaches $(g, k)$ as in Theorem \ref{thm:density}, bearing in mind Proposition \ref{prop:densitytrace}. The results in Section \ref{sec:conformalstructure} apply to $(g^j, k^j)$. We can choose $\Lambda_\beta, c_\beta$, and $b_\beta$ (with $\beta = 1 + q$ when $n >3$ and $\beta \in (2, \min \{3, \gamma\})$ when $n=3$) in Proposition \ref{prop:barriers} so that the conclusions there hold for all sufficiently large $j$. The graphs of the functions $f_{\Sigma^j}$ from Proposition \ref{prop:limitanalysis} (\ref{exteriorsolution}) are $2 C$-minimizing boundaries in $(M \times \R, g^j + d t^2)$ where $C := 1 + n \sup_j \sup_{x \in M} |k_j(x)|_{g^j}$. The $W^{3, p}_{1-q}$ estimates for $f_{\Sigma^j}$ on $\{x \in M : |x| > 2 \Lambda_\beta\}$ are uniform in $j$. The $\C^{3, \alpha}_{loc}$ estimates are locally uniform in $j$.  \\

Passing to a geometric limit along a subsequence of $\{j\}$ (which we continue to denote by $\{j\}$), we conclude that there exist $\Sigma \subset M \times \R$, $U_\Sigma \subset M$, $f_\Sigma : U_\Sigma \to \R$ as in Proposition \ref{prop:limitanalysis} (\ref{exteriorsolution}) such that $\Sigma^j \to \Sigma$ as single-layered $\C^{3, \alpha}_{loc}$ hypersurfaces of $M \times \R$ and such that $\partial U_{\Sigma^j} \to \partial U_\Sigma$ as single-layered $\C^{3, \alpha}_{loc}$ hypersurfaces of $M$. \\

From Proposition \ref{prop:SchoenYauidentity} we know that the Yamabe type of the components of $\partial U_{\Sigma^j}$ is positive, because the strict dominant energy condition holds for $(g^j, k^j)$. The Yamabe type of the components of $\partial U_\Sigma$ might be zero though.\\

Let $t_0^j \nearrow \infty$ be a sequence such that $\pm t_0^j$ are regular values for both $f_\Sigma$ and $f_{\Sigma^j}$ for every $j$. Let $\tilde g^j$ be metrics on $\Sigma^j$ as in Proposition \ref{prop:normalizedSigma} such that $\tilde g^j = \bar g^j$ on $\Sigma^j \cap (M \times (-t_0^j, t_0^j))$. Let $u^j \in \C^{2, \alpha}_{loc} (\Sigma^j)$ be the solutions of $- \Delta_{\tilde g^j} u^j+ c_n R_{\tilde g^j} u^j = 0$ from Proposition \ref{prop:conformalfactor}. Let $E^j$ denote the energy of $(g^j, k^j)$, and let $a^j_1 \leq 0$ be as in Proposition \ref{prop:conformalfactor}. From (\ref{eqn:a1}) we see that $\int_{N^j} |d u^j|_{\bar g^j}^2 \leq 2 (2-n) a_1^j |\mathbb{S}^{n-1}|$. By assumption, $E^j \to 0$ and hence $a^j \to 0$ as $j \to \infty$. The proof of Proposition \ref{prop:main} gives that $a^j + a^j_1 \geq 0$. It follows that $a^j_1 \to 0$. In conjunction with the Sobolev inequality and the equation that $u^j$ satisfies we see that $u^j (x) \to 1$ as $|x| \to \infty$ uniformly in $j$. Using standard elliptic theory we conclude that $u_j$ converges in $\C^{2, \alpha}_{loc}$ to the constant function one on $\Sigma$. Applying (\ref{eqn:a1}) once more we deduce that $R_{\bar g} = 0$ and that $h = k$ on $\Sigma$. \\ 

Using our decay assumptions for $f_\Sigma$ and that $R_\Sigma = 0$, we see that the metric $\bar g = g + df_\Sigma \otimes df_\Sigma$ on $\Sigma$ is asymptotically flat of type $(2, p, q, q_0, \alpha)$ on $\{x \in N : |x| > \Lambda_\beta\}$ for every $q_0>0$ and that the energy of the end $N$ with respect to $\bar g$ vanishes. Here we view $(\Sigma, \bar g)$ as a time-symmetric initial data set, i.e. one with zero spacetime second fundamental form. We introduce a new distance function $0 < s \in \C^{3, \alpha}_{loc} (\Sigma)$ that agrees with $|x|$ on $\{x \in M : |x| > 2 \Lambda_\beta\}$ and such that $s(x) = |t|^{-1}$ on $\Sigma \cap \{ (x, t) : x \in M \}$ for $|t|$ large. \\

Next, we show that $\Ric_{\bar g} \equiv 0$. We follow the ideas of R. Schoen and S.-T. Yau closely (see in particular \cite[p. 72-74]{Schoen-Yau:1979-pmt1} and \cite[Lemma 3.2]{Schoen-Yau:1979-pmt1}), making technical adjustments to accommodate that $(\Sigma, \bar g)$ may contain cylindrical ends and that we assume less regularity for $\bar g$ than in \cite{Schoen-Yau:1979-pmt1}. Let $h \in \C^{2, \alpha}_{c} (\text{Sym}^2(T^*\Sigma))$ be a compactly supported symmetric $(0, 2)$-tensor. For small values of $\kappa$, consider the metric $\bar g_\kappa = \bar g + \kappa h$. Let $\sigma_0 >0$ small be such that for all $\sigma \in (0, \sigma_0)$, both $\sigma$ and $\sigma^{-1}$ are regular values of $s(x)$.  Let $0 \leq q \in \C^{2, \alpha}(\Sigma)$ be a function that coincides with $|x|^{{-2n}}$ on $\{ |x| > 2 \Lambda_\beta\}$, and such that $\supp (q) \cap \{s(x) < \sigma_0\} = \emptyset$. For $\sigma \in (0, \sigma_0)$ and sufficiently small $\kappa$, we can solve the mixed Dirichlet\slash Neumann problems 
\begin{eqnarray*}
 \left\{ \begin{array}{lllll} - \Delta_{\bar g_\kappa} u_{\kappa, \sigma} + c_n R_{\bar g_\kappa} u_{\kappa,\sigma} &=& \kappa^2 q &\text{ on }& \{\sigma < s(x) < \sigma^{-1}\} \\ u_{\kappa, \sigma} &=& 1 &\text{ on }& \{s(x) = \sigma^{-1}\} \\ \nu_{\bar g_\kappa} (u_{\kappa, \sigma}) &=& 0 &\text{ on }& \{s(x) = \sigma \}. \end{array} \right.
\end{eqnarray*}
To see this, note that if $w$ is a solution of the corresponding homogeneous problem (with zero right-hand side), then we can multiply the equation satisfied by $w$ on $\{\sigma < s(x) < \sigma^{-1} \}$, integrate by parts, and use the Sobolev inequality on $(\{x \in \Sigma : s(x) > \sigma_0\}, \bar g_\kappa)$ to obtain that 
\begin{eqnarray*}
&& \left( \int_{\{ \sigma_0 < s(x) < \sigma^{-1}\}} |w|^{\frac{2n}{n-2}} d \L^n_{\bar g_\kappa}  \right)^{\frac{n-2}{n}}   \leq C \int_{\{ \sigma_0 < s(x) < \sigma^{-1}\}} |d w|_{\bar g_\kappa}^2 d \L^n_{\bar g_\kappa} \\ &=&  C \int_{\{ \sigma_0 < s(x) < \sigma^{-1}\}} - R_{\bar g_\kappa} w^2 d \L^n_{\bar g_\kappa} \\ &\leq& C \left( \int_{\{ \sigma_0 < s(x) < \sigma^{-1}\}} |R_{\bar g_\kappa}|^{\frac{n}{2}} d \L^n_{\bar g_\kappa}  \right)^{\frac{2}{n}} \left( \int_{\{ \sigma_0 < s(x) < \sigma^{-1}\}} |w|^{\frac{2n}{n-2}} d \L^n_{\bar g_\kappa}\right)^{\frac{n-2}{n}}. 
\end{eqnarray*}
Here, $C$ depends only on a uniform upper bound for the Sobolev constant of $(\{x \in \Sigma : s(x) > \sigma_0\}, \bar g_\kappa)$ for all $|\kappa|$ small. Clearly, $||R_{\bar g_\kappa}||_{L^{\frac{n}{2}}} = O(|\kappa|)$, so that $w = 0$ on $\{x \in \Sigma : s(x) > \sigma_0\}$ when $|\kappa|$ is small. Using the equation that $w$ satisfies, we conclude that $w$ vanishes on all of $\{\sigma < s(x) < \sigma^{-1}\}$, provided again that $|\kappa|$ is sufficiently small. \\ 

Let $v_{\kappa, \sigma} := u_{\kappa, \sigma} -1$. A similar argument as above shows that 
\begin{eqnarray*}
&& \left( \int_{\{ \sigma_0 < s(x) < \sigma^{-1}\}} |v_{\kappa, \sigma}|^{\frac{2n}{n-2}} d \L^n_{\bar g_\kappa}  \right)^{\frac{n-2}{n}}  \leq C \int_{\{ \sigma < s(x) < \sigma^{-1}\}} |d v_{\kappa, \sigma}|_{\bar g_\kappa}^2 d \L^n_{\bar g_\kappa} \\  =  && C  \int_{\{ \sigma < s(x) < \sigma^{-1}\}} - c_n R_{\bar g_\kappa} v_{\kappa, \sigma}^2 + (\kappa^2 q - c_n R_{\bar g_\kappa}) v_{\kappa, \sigma} d \L^n_{\bar g_\kappa} \\ 
\leq &&  C \left( \int_{ \{\sigma_0 < s(x) < \sigma\}} |R_{\bar g_\kappa}|^{\frac{n}{2}} d \L^n_{\bar g_\kappa} \right) ^{\frac{n}{2}} \left( \int_{ \{\sigma_0 < s(x) < \sigma\}} |v_{\kappa, \sigma}|^{\frac{2n}{n-2}} d \L^n_{\bar g_\kappa} \right)^{\frac{n-2}{n}} \\ && + C \left( \int_{ \{\sigma_0 < s(x) < \sigma\}} |\kappa^2 q - c_n R_{\bar g_\kappa}|^{\frac{2n}{n+2}} d \L^n_{\bar g_\kappa} \right)^{\frac{n+2}{2n}} \left( \int_{\{\sigma_0 < s(x) < \sigma \}} |v_{\kappa, \sigma}|^{\frac{2n}{n-2}} d \L^n_{\bar g_\kappa}\right)^{\frac{n-2}{2n}}. 
\end{eqnarray*}
This estimate implies that $||v_{\kappa, \sigma}||_{L^{2n/(n-2)} (\{\sigma_0 < s(x) < \sigma^{-1}\})} = O(|\kappa|)$. Reasoning as in the proof of Proposition \ref{prop:conformalfactor}, we obtain that $|| v_{\kappa, \sigma} ||_{\C^{2, \alpha} (\{2 \sigma < s(x) < (2 \sigma)^{-1}\})} = O(|\kappa|)$, where the quantity on the right is independent of $\sigma \in (0, \sigma_0)$. \\

Let $u_\kappa$ be a subsequential limit  of $u_{\kappa, \sigma}$ as $\sigma \searrow 0$. Then $||u_\kappa - 1||_{\C^{2, \alpha} (\Sigma)} = O(|\kappa|)$, $- \Delta_{\bar g_\kappa} u_\kappa + c_n R_{\bar g_\kappa} u_\kappa = \kappa^2 q$, and $u(x) \to 1$ as $|x| \to \infty$. Since each $u_{\kappa, \sigma}$ is harmonic on $\{\sigma < s(x) < \sigma_0\}$ and satisfies a Neumann boundary condition on $\{s(x) = \sigma\}$, it follows that $\int_{\{s(x) = \sigma_0\}} \nu_{\bar g_\kappa} (u) d \H^{n-1}_{\bar g_\kappa} = 0$. Asymptotic analysis as in e.g. \cite{Bartnik:1986} shows that $u_\kappa(x) = 1 + A_\kappa |x|^{2-n} + O^{2, \alpha} (|x|^{2-n - \eta})$ as $|x| \to \infty$ for constants $A_\kappa \in \R$ and any $\eta \in (0, q)$. An integration by parts shows that $(n-2) |\mathbb{S}^{n-1}| A_\kappa = \int_\Sigma (\kappa^2 q - c_n R_{\bar g_\kappa} u_\kappa) d \L^n_{\bar g_\kappa}$.  As in \cite{Schoen-Yau:1979-pmt1}, using that  $\bar g_0 = \bar g$, that $R_{\bar g_0} = 0$, that $u_0 = 1$, and that $||u_\kappa - 1||_{\C^{2, \alpha} (\Sigma)} = O(|\kappa|)$, we see that $A_\kappa$ is differentiable at $\kappa = 0$, and that its derivative equals
\begin{eqnarray} \label{eqn:ddta} 4 (n-1) |\mathbb{S}^{n-1}| \frac{d}{d\kappa} |_{\kappa=0} A_\kappa &=& - \int_\Sigma \frac{d}{d\kappa}|_{\kappa=0} R_{\bar g_\kappa} d \L^n_{\bar g} \\ &=& \int_\Sigma \Delta_{\bar g} \text{ tr}_{\bar g} (h) - \div_{\bar g } \div_{\bar g} h  + \bar g (h, \Ric_{\bar g}) d \L^n_{\bar g} \nonumber \\ &=& \int_{\Sigma} \bar g (h, \Ric_{\bar g}) d \L^n_{\bar g}. \nonumber
\end{eqnarray} 
The scalar curvature of $u_\kappa^{4/(n-2)} \bar g_\kappa$ is non-negative, and positive on $\{ |x| > 2 \Lambda_\beta \}$ when $\kappa \neq 0$. Using the same additional argument as in the proof of Proposition \ref{prop:main} we can justify the use of the Riemannian positive energy theorem.\footnote{There is an additional subtlety here. The construction of barriers in the end for minimizing hypersurfaces in \cite{Schoen-Yau:1979-pmt1, Schoen:1989} requires that the metric has harmonic asymptotics. The perturbation to harmonic asymptotics requires an initial step that is described in detail in \cite{Schoen-Yau:1981-asymptotics}. It proceeds by a linear homotopy of the metric to the Euclidean metric far out in the end and a conformal perturbation to reimpose the constraint, all while changing the mass by no more than some prescribed amount. This step can be integrated into the preceding perturbation argument. The necessary modifications are minor. The decay of $u_\kappa$ is then improved to $1 + A_\kappa |x|^{2-n} + O^{2, \alpha}(|x|^{1-n})$.} Thus the energy of the end $N$ with respect to the metrics $u_\kappa^{4/(n-2)} \bar g_\kappa$ is non-negative. In view of the expansion of $u_\kappa(x)$ as $|x| \to \infty$ and the fact that the energy of $\bar g$ vanishes, we see that the energy of the asymptotically flat end of the metric $u_\kappa^{4/(n-2)} \bar g_\kappa$ equals $\frac{n-2}{2} A_\kappa$. It follows that $\frac{d}{d\kappa}|_{\kappa=0} A_\kappa = 0$ and hence, from (\ref{eqn:ddta}), that $\int_{\Sigma} \bar g( h, \Ric_{\bar g}) d \L^n_{\bar g} = 0$. \\

Let $\chi \in \C_c^{3, \alpha} (\Sigma)$ be a non-negative function with support in a single coordinate chart. Let $h_i \in \C^{2, \alpha}_c (\text{Sym}^2 (T^*\Sigma))$ be a sequence of symmetric $(0, 2)$-tensors that approximate $\chi \Ric_{\bar g}$ in $\C^{0, \alpha} (\Sigma)$. We know that $\int_{\Sigma} \bar g (h_i, \Ric_{\bar g}) d \L^n_{\bar g} = 0$ for all $i = 1, 2, \ldots$. Passing to the limit as $i \to \infty$, we see that $\int_{\Sigma} \chi |\Ric_{\bar g}|_{\bar g}^2 d \L^n_{\bar g} = 0$. Thus $\Ric_{\bar g} = 0$. In particular, $(\Sigma, \bar g)$ is analytic \cite[Theorem 5.26]{Besse:2008}.  \\
 
If $\Sigma$ has cylindrical ends, then it is an easy matter to construct a geodesic line in $\Sigma$. By the Cheeger-Gromoll splitting theorem, $\Sigma$ splits off a factor of $\R$ isometrically. This is clearly impossible. Thus $\Sigma$ has no cylindrical ends. By the Bishop-Gromov comparison theorem, for any $x \in \Sigma$, the function $r \to  (\omega_n r^n)^{-1}  \L^n_{\bar g}(B_{\bar g}(x, r))$ is non-increasing. Explicit comparison with coordinate balls in the asymptotically flat end of $\Sigma$ shows that as $r \to \infty$, this quantity converges to $1$. It follows that $\L^n_{\bar g}(B_{\bar g}(x, r)) = \omega_{n} r^n $ for all $r\geq0$. Equality in the Bishop-Gromov comparison theorem holds only for geodesic balls in the model space. Thus $(\Sigma, \bar g)$ is isometric to Euclidean space $(\R^n, \sum_{i=1}^n dx_i^2)$. \\

The last step of the argument is exactly as in \cite[p. 260]{Schoen-Yau:1981-pmt2}. (Reading this step backwards serves as motivation for introducing the Jang equation in the first place, cf. \cite{Jang:1978}.) We think of $f_\Sigma$ as a function on $\Sigma$. Let $(y_1, \ldots, y_n)$ be a Euclidean coordinate system on $\Sigma$. Then $\bar g_{ij} = \delta_{ij}$ and $g_{ij} = \delta_{ij} - (f_\Sigma)_i (f_\Sigma)_j$. Recall that $h_{ij} = k_{ij}$. A computation shows that $h_{ij} =  (1 + |df_\Sigma|^2_{g})^{-1/2} (\nabla_{g}^2 f_\Sigma)_{ij} = (1 + |df_\Sigma|^2_{g})^{-1/2} \partial^2_{ij} f_\Sigma $. It follows that $g$ and $k$ agree, respectively, with the pull-back metric and the second fundamental form of the embedding $M \to (\R^n \times \R,  dx_1^2 + \ldots + dx_n^2 - dx_{n+1}^2)$ defined by $(y_1, \ldots, y_n) \to (y_1, \ldots, y_n, f_\Sigma(y_1, \ldots, y_n))$.          
\end{proof}
\end{proposition}

\begin{remark} M. Nardmann \cite{Nardmann:2010} has shown that initial data sets that satisfy the Gauss and Codazzi equations for constant curvature spaces may be imbedded therein. In particular, initial data sets $(M, g, k)$ with $\mu = 0$, $J = 0$ are submanifolds of Minkowski space. 
\end{remark}


\appendix 

\section{}

The point of the following version of the Sobolev inequality is that no particular boundary behavior is assumed for $\phi$ on $\partial M$. 

\begin{lemma} [\protect{Cf. \cite[Lemma 3.1]{Schoen-Yau:1979-pmt1}}] 
\label{lem:SobolevInequality}
Let $(M, g)$ be a complete connected Riemannian manifold, possibly with boundary, such that there exists a compact set $K \subset M$ and a diffeomorphism $x = (x_1, \ldots, x_n) : M \setminus K \to \R^n \setminus \bar B(0, 1)$ so that for some constant $c \geq 1$ we have that $c^{-1} \delta_{ij} \leq g_{ij}  \leq c \delta_{ij}$, as quadratic forms. For every $1 \leq p < n$ there exists a constant $C = C(M, g, p)$ such that 
  \begin{equation*}
    \left( \int_{ M}
      |\phi|^{\frac{np}{n-p}} d \L^n_{ g}\right)^{\frac{n-p}{np}} \leq C
   \left(  \int_{ M} |d \phi|_g^p d \L^n_{ g} \right)^{\frac{1}{p}} \text{ for all } \phi \in
    \C_c^1(M).
  \end{equation*}
 \end{lemma}
 
\begin{lemma} \label{lem:monotonicity} Let $r>0$ and $c \geq 1$ be two constants and let $g_{ij}$ be the components of a metric on $B(0, 2r) \subset \R^n$ such that $c^{-1} \delta_{ij} \leq g_{ij} \leq c \delta_{ij}$ as quadratic forms. If $T$ is a $g$--area minimizing boundaryless $(n-1)$-dimensional integer multiplicity current in $B(0, 2r)$ with $0 \in \supp (T)$, then $M_{B(0, r)}^{g} (T)\geq c^{- 3 (n-1)^2} \omega_{n-1} r^{n-1}$ where the mass is computed with respect to the $g$--metric. \begin{proof} The current $T$ satisfies the conditions of Lemma 5.1 in \cite{Bray-Lee:2009} with $\gamma = c^{2(n-1)}$. \end{proof}
\end{lemma}

\bibliographystyle{amsplain}
\bibliography{JangReductionreferences}

\providecommand{\bysame}{\leavevmode\hbox to3em{\hrulefill}\thinspace}
\providecommand{\MR}{\relax\ifhmode\unskip\space\fi MR }
\providecommand{\MRhref}[2]{%
  \href{http://www.ams.org/mathscinet-getitem?mr=#1}{#2}
}
\providecommand{\href}[2]{#2}
\begin{thebibliography}{10}

\bibitem{Almgren:1976}
F.~J. Almgren, Jr., \emph{Existence and regularity almost everywhere of
  solutions to elliptic variational problems with constraints}, Mem. Amer.
  Math. Soc. \textbf{4} (1976), no.~165, viii+199. \MR{0420406 (54 \#8420)}

\bibitem{Andersson-Eichmair-Metzger:2010}
Lars Andersson, Michael Eichmair, and Jan Metzger, \emph{Jang's equation and
  its applications to marginally trapped surfaces}, Complex analysis and
  dynamical systems {IV}. {P}art 2, Contemp. Math., vol. 554, Amer. Math. Soc.,
  Providence, RI, 2011, pp.~13--45. \MR{2884392}

\bibitem{Andersson-Metzger:2009}
Lars Andersson and Jan Metzger, \emph{The area of horizons and the trapped
  region}, Comm. Math. Phys. \textbf{290} (2009), no.~3, 941--972. \MR{2525646
  (2010f:53118)}

\bibitem{ADM:1962}
R.~Arnowitt, S.~Deser, and C.~W. Misner, \emph{The dynamics of general
  relativity}, Gravitation: {A}n introduction to current research, Wiley, New
  York, 1962, pp.~227--265. \MR{0143629 (26 \#1182)}

\bibitem{ADM:2004}
\bysame, \emph{The {D}ynamics of {G}eneral {R}elativity}, arXiv:gr-qc/0405109v1
  (2004).

\bibitem{Bartnik:1986}
Robert Bartnik, \emph{The mass of an asymptotically flat manifold}, Comm. Pure
  Appl. Math. \textbf{39} (1986), no.~5, 661--693. \MR{849427 (88b:58144)}

\bibitem{Besse:2008}
Arthur~L. Besse, \emph{Einstein manifolds}, Classics in Mathematics,
  Springer-Verlag, Berlin, 2008, Reprint of the 1987 edition. \MR{2371700
  (2008k:53084)}

\bibitem{Bray-Lee:2009}
Hubert~L. Bray and Dan~A. Lee, \emph{On the {R}iemannian {P}enrose inequality
  in dimensions less than eight}, Duke Math. J. \textbf{148} (2009), no.~1,
  81--106. \MR{2515101 (2010f:53051)}

\bibitem{CO:1981}
D.~Christodoulou and N.~{\'O}~Murchadha, \emph{The boost problem in general
  relativity}, Comm. Math. Phys. \textbf{80} (1981), no.~2, 271--300.
  \MR{623161 (84e:83011)}

\bibitem{Deser:1983}
S.~Deser, \emph{Positive classical gravitational energy from classical
  supergravity}, Phys. Rev. D (3) \textbf{27} (1983), no.~12, 2805--2808.
  \MR{708019 (84m:83063)}

\bibitem{Deser-Teitelboim:1977}
S.~Deser and C.~Teitelboim, \emph{Supergravity has positive energy}, Phys. Rev.
  Lett. \textbf{39} (1977), no.~5.

\bibitem{Duzaar-Steffen:1993}
Frank Duzaar and Klaus Steffen, \emph{{$\lambda$} minimizing currents},
  Manuscripta Math. \textbf{80} (1993), no.~4, 403--447. \MR{1243155
  (95f:49062)}

\bibitem{Eichmair:2009-Plateau}
Michael Eichmair, \emph{The {P}lateau problem for marginally outer trapped
  surfaces}, J. Differential Geom. \textbf{83} (2009), no.~3, 551--583.
  \MR{2581357}

\bibitem{Eichmair:2010}
\bysame, \emph{Existence, regularity, and properties of generalized apparent
  horizons}, Comm. Math. Phys. \textbf{294} (2010), no.~3, 745--760.
  \MR{2585986 (2011d:53171)}

\bibitem{spacetimePMT}
Michael Eichmair, Lan-Hsuang Huang, Dan Lee, and Richard Schoen, \emph{The
  spacetime positive mass theorem in dimensions less than eight},
  arXiv:1110.2087.

\bibitem{potpourri}
Michael Eichmair and Jan Metzger, \emph{Jenkins-{S}errin type results for the
  {J}ang equation}, arXiv:1205.4301.

\bibitem{Galloway-Schoen:2006}
Gregory~J. Galloway and Richard Schoen, \emph{A generalization of {H}awking's
  black hole topology theorem to higher dimensions}, Comm. Math. Phys.
  \textbf{266} (2006), no.~2, 571--576. \MR{2238889 (2007i:53078)}

\bibitem{Gilbarg-Trudinger:1998}
David Gilbarg and Neil~S. Trudinger, \emph{Elliptic partial differential
  equations of second order}, Classics in Mathematics, Springer-Verlag, Berlin,
  2001, Reprint of the 1998 edition. \MR{1814364 (2001k:35004)}

\bibitem{Grisaru:1978}
M.~T. Grisaru, \emph{Positivity of the energy in {E}instein theory}, Phys.
  Lett. \textbf{73B} (1978), no.~2.

\bibitem{Jang:1978}
Pong~Soo Jang, \emph{On the positivity of energy in general relativity}, J.
  Math. Phys. \textbf{19} (1978), no.~5, 1152--1155. \MR{488515 (80b:83012a)}

\bibitem{Metzger:2010}
Jan Metzger, \emph{Blowup of {J}ang's equation at outermost marginally trapped
  surfaces}, Comm. Math. Phys. \textbf{294} (2010), no.~1, 61--72. \MR{2575475
  (2011a:53139)}

\bibitem{Nardmann:2010}
Marc Nardmann, \emph{A remark on the rigidity case of the positive energy
  theorem}, arXiv:1004.5430v1 [math.DG] (2010).

\bibitem{Schoen-Simon-Yau:1975}
R.~Schoen, L.~Simon, and S.~T. Yau, \emph{Curvature estimates for minimal
  hypersurfaces}, Acta Math. \textbf{134} (1975), no.~3-4, 275--288.
  \MR{0423263 (54 \#11243)}

\bibitem{Schoen:1984}
Richard Schoen, \emph{Conformal deformation of a {R}iemannian metric to
  constant scalar curvature}, J. Differential Geom. \textbf{20} (1984), no.~2,
  479--495. \MR{788292 (86i:58137)}

\bibitem{Schoen:1989}
\bysame, \emph{Variational theory for the total scalar curvature functional for
  {R}iemannian metrics and related topics}, Topics in calculus of variations
  ({M}ontecatini {T}erme, 1987), Lecture Notes in Math., vol. 1365, Springer,
  Berlin, 1989, pp.~120--154. \MR{994021 (90g:58023)}

\bibitem{Schoen-Simon:1981}
Richard Schoen and Leon Simon, \emph{Regularity of stable minimal
  hypersurfaces}, Comm. Pure Appl. Math. \textbf{34} (1981), no.~6, 741--797.
  \MR{634285 (82k:49054)}

\bibitem{Schoen-Yau:1979-pmt1}
Richard Schoen and Shing~Tung Yau, \emph{On the proof of the positive mass
  conjecture in general relativity}, Comm. Math. Phys. \textbf{65} (1979),
  no.~1, 45--76. \MR{526976 (80j:83024)}

\bibitem{Schoen-Yau:1981-asymptotics}
\bysame, \emph{The energy and the linear momentum of space-times in general
  relativity}, Comm. Math. Phys. \textbf{79} (1981), no.~1, 47--51. \MR{609227
  (82j:83045)}

\bibitem{Schoen-Yau:1981-pmt2}
\bysame, \emph{Proof of the positive mass theorem. {II}}, Comm. Math. Phys.
  \textbf{79} (1981), no.~2, 231--260. \MR{612249 (83i:83045)}

\bibitem{Schoen-Yau:1979-pat}
Richard~M. Schoen and Shing~Tung Yau, \emph{Complete manifolds with nonnegative
  scalar curvature and the positive action conjecture in general relativity},
  Proc. Nat. Acad. Sci. U.S.A. \textbf{76} (1979), no.~3, 1024--1025.
  \MR{524327 (80k:58034)}

\bibitem{Tamanini:1982}
Italo Tamanini, \emph{Boundaries of {C}accioppoli sets with
  {H}\"older-continuous normal vector}, J. Reine Angew. Math. \textbf{334}
  (1982), 27--39. \MR{667448 (83m:49067)}

\bibitem{Witten:1981}
Edward Witten, \emph{A new proof of the positive energy theorem}, Comm. Math.
  Phys. \textbf{80} (1981), no.~3, 381--402. \MR{626707 (83e:83035)}

\end{thebibliography}
\end{document}